\documentclass[a4paper,fleqn]{cas-sc}

\usepackage{amssymb}
\usepackage{amsmath}
\usepackage{amsthm}
\usepackage{amsfonts}
\usepackage{graphicx}
\usepackage{color}
\usepackage{subcaption}
\usepackage[T1]{fontenc}
\usepackage{algorithm}
\usepackage{algpseudocode}
\usepackage[numbers]{natbib}

\makeatletter
\AtEndDocument{\clearpage\addtocounter{page}{-1}\immediate\write\@auxout{\string\csxdef{lastpage}{\thepage}}\addtocounter{page}{1}}
\makeatother

\theoremstyle{remark}
\newtheorem{remark}{Remark}
\numberwithin{remark}{subsection}

\begin{document}

\shorttitle{Remap for ALE Hydrodynamics}
\shortauthors{Keith et al.}



\title{Remap for ALE Hydrodynamics via
		   Direct Interpolation and Constraint Optimization}

\author[bru]{Brendan Keith}
\author[bru]{Dohyun Kim}
\author[lab]{Boyan Lazarov}
\author[lab]{Ketan Mittal}
\author[lab]{Cosmin G. Petra}
\author[lab]{Mathias R. Schmidt}
\author[lab]{Vladimir Z. Tomov}
\cormark[1]
\cortext[1]{Corresponding author: tomov2@llnl.gov}

\affiliation[lab]{organization={Lawrence Livermore National Laboratory},
	city={Livermore}, state={CA}, country={USA}}
\affiliation[bru]{organization={Brown University},
	city={Providence}, state={RI}, country={USA}}

\begin{abstract}
This work presents an optimization-based remap approach for arbitrary
Lagrangian–Eulerian (ALE) hydrodynamics that preserves
sharp features and physical bounds while conserving global quantities such as
volume, mass, momentum, and total energy.
Achieving these properties is the main challenge in remap methods.
The presented approach is well-suited for high-order finite element (FE)
discretizations and ALE methods that allow for large, purely
Lagrangian displacements between successive remeshing and remapping steps.
Recent advances in parallel, GPU-capable interpolation routines
enable the initial transfer of all fields via direct physical-space interpolation,
avoiding computationally costly and diffusive pseudo-transient PDE-based advection.
The interpolation produces bounded fields with minimal numerical diffusion.
However, it may violate the desired conservation of global properties.
Thus, a subsequent optimization step restores
conservation while preserving physical bounds.
The presented method provides minimally diffusive remapping of
finite element fields of arbitrary order.
Moreover, it directly remaps quantities defined at integration points,
eliminating the need for projection onto finite element spaces.
The approach also mitigates typical challenges associated with synchronizing
primal variables and conserved quantities in advection-based remap methods.
The proposed optimization-based remap method is evaluated using standard two-
and three-dimensional remap benchmarks and full ALE hydrodynamics simulations.
\end{abstract}



\begin{keywords}
		ALE remap \sep finite elements \sep high order \sep optimization



\end{keywords}

\maketitle



\section{Introduction}
\label{sec_intro}

This work presents a novel optimization-based remap approach for
operator-split arbitrary Lagrangian-Eulerian (ALE) hydrodynamics.
An operator-split ALE simulation generally consists of three stages:
(i) Lagrangian evolution, in which the hydrodynamic equations are
advanced in time while the mesh moves with the material;
(ii) remeshing, in which the mesh is modified to improve its quality; and
(iii) remapping, in which all solution fields are transferred from
the deformed Lagrangian mesh to the modified mesh.
This manuscript focuses exclusively on the remapping process.
The remap phase of an ALE method must produce accurate and sharp fields with
minimal numerical diffusion around steep solution gradients,
conserve physical quantities such as mass and energy,
and preserve the local minimum and maximum bounds of the remapped fields.
Moreover, in the context of this work's starting discretization method
\cite{Tomov2018}, the remap algorithm must support arbitrary-order
finite element (FE) discretizations and curved meshes.

Most existing remap methods in ALE frameworks transfer solution fields between
meshes through pseudo-transient advection-based processes, where flux-based
quantities ensure conservation.
Local bounds are enforced using additional nonlinear techniques such as
Flux-Corrected Transport \cite{Dobrev2015, Manuel2017, Tomov2018} or
optimization-based formulations introduced by \cite{Liska2010} and later extended
by Bochev et al.~\cite{Bochev2011, Bochev2013, Bochev2014, Bochev2020, Kuzmin2024}.
These advection-based approaches inherently introduce numerical diffusion that
accumulates over successive pseudotime steps, progressively degrading sharp
solution gradients.
This diffusion is the primary source of numerical error in existing ALE methods,
compromising the key advantage of Lagrangian hydrodynamics:
the preservation of sharp material interfaces.

This work explores a more direct remapping approach that
eliminates the need for pseudotime time-stepping.
Given the initial and optimized meshes,
direct interpolation is first applied in physical space
using the capabilities of the recently developed general field evaluation framework \cite{Mittal2025}.
This step does not involve the construction of numerical fluxes; hence, it does not
require the two meshes to be close to each other or share the same topology.
The interpolated field has the maximum possible sharpness while preserving
local bounds, though it initially violates conservation.
To address this, a second step applies an optimization-based correction that
restores conservation while maintaining local bounds and minimizing
the deviation from the interpolated solution to preserve its sharpness.
Two alternative optimization strategies are considered:
an interior-point method (IPM) implemented in the HiOp library \cite{Petra2019}
and an entropic Bregman projection algorithm inspired by recent work
incorporating entropic regularization into computational mechanics
\cite{keith2026proximal,dokken2025latent,kim2025simple}.
Furthermore, the proposed approach enables direct remapping of point-based quantities,
such as fields defined only at quadrature points, without conversion to and from a
finite element representation, avoiding additional diffusion and numerical errors.

The paper extends the scalar-field remap framework introduced in \cite{Remap2026}
by focusing on single-material settings where the state is described by an
indicator field $\eta \in [0,1]$, along with the associated hydrodynamic fields.
The indicator field represents the material support and accurately describes sharp
transitions that are particularly sensitive to numerical diffusion during remap.
This setting allows to assess the coupling among solution fields and evaluate
the accuracy and sharpness of the proposed interpolation-optimization workflow.
It also avoids the additional constraints associated with multiple materials.
Although the proposed remapping approach is presented
in the context of ALE hydrodynamics, it is applicable to other
Lagrangian simulations that require field transfer between meshes,
including simulations of fluid flow, elastoplastic deformation,
or molecular dynamics \cite{Scrimieri2014, Bussetta2015, Levy2024}.

The remainder of this paper is organized as follows.
The problem formulation is presented in Section~\ref{sec_problem}.
Section~\ref{sec_interp} summarizes the basic concepts of the interpolation step.
Section~\ref{sec_bounds} discusses the definition of the admissible solution bounds.
Details of the optimization solvers are provided in Section~\ref{sec_optimization}.
Numerical examples are provided in Section~\ref{sec_results}.
Section~\ref{sec_concl} summarizes the work and presents final conclusions.


\section{Problem formulation}
\label{sec_problem}

Remap is considered within an ALE workflow, where solution fields are
transferred from an initial mesh $\Omega_h^0$ to a new mesh $\Omega_h$
that typically results from a mesh optimization step.
The focus is a single-material state described by a material indicator $\eta$,
density $\rho$, specific internal energy $e$, and velocity $v$.
The indicator function $\eta$ represents the material support and
can contain sharp transitions, thus
controlling numerical diffusion and preserving its local bounds is essential.

The initial solution state is denoted by $(\eta^0, \rho^0, e^0, v^0)$,
associated with the initial mesh $\Omega_h^0$.
The state fields follow the discretization choices of the underlying Lagrangian
hydrodynamics framework~\cite{Dobrev2012,Laghos}:
$\eta^0$, $\rho^0$, and $e^0$ are represented as quadrature point values,
while $v^0$ is represented as a vector-valued, $H^1$-conforming FE function.
Throughout this manuscript, fields represented by values at quadrature points
are referred to as Q-functions.

Although pressure is a derived thermodynamic quantity,
controlling it during remap is important for robustness:
pressure enters the hydrodynamic update directly, and nonphysical pressure
values can destabilize the subsequent Lagrangian evolution.
In particular, bounds on specific internal energy $e$ alone
do not imply pointwise bounds on pressure, because pressure depends on both
density and internal energy through the nonlinear equation of state.
Pressure is therefore constructed as an auxiliary Q-function,
$p^0 := p(\rho^0,e^0)$, from the input Q-functions using the
hydrodynamic equation of state, and remap $p$ alongside $\eta$ and $\rho$.
The optimization then enforces the pressure bounds directly;
after remap, the specific internal energy is recovered as $e=e(\rho,p)$,
returning the state in its original variables.

Given the initial and optimized meshes, the approach consists of two steps:
(i) direct interpolation in physical space,
which produces sharp, locally bounded fields; and
(ii) an optimization-based correction that
restores global conservation while maintaining local bounds.
The interpolation procedure transfers all state fields from the initial mesh
$\Omega_h^0$ to the optimized mesh $\Omega_h$,
yielding the interpolated states  $(\eta^*, \rho^*, p^*, v^*)$
as detailed in Section~\ref{sec_interp}.
These states serve as both the initial guess and the target
for the subsequent optimization-based reconstruction.
The optimization problems minimize the sum of the $L^2$-norms of
$\eta-\eta^*$, $\rho-\rho^*$, and $p-p^*$ and the $H_{1}$-norm of $v-v^*$,
subject to conservation of the volume, mass, momentum, and total energy:

\begin{subequations}
\label{eq_problem}
\begin{align}
\label{eq_obj}
  & \min_{\eta, \rho, p,v} ~ ||\eta - \eta^*||_{L^2}^2 +
                           ||\rho - \rho^*||_{L^2}^2 +
                           ||p - p^*||_{L^2}^2+
                           ||v - v^*||_{H_1}^2  \\
  & \quad\text{subject to:} \nonumber  \\
\label{eq_vol}
  & \int_{\Omega_h} \eta = \int_{\Omega_h^0} \eta^0,      \\
\label{eq_mass}
  & \int_{\Omega_h} \eta \rho = \int_{\Omega_h^0} \eta^0 \rho^0,  \\
\label{eq_momentum}
  & \int_{\Omega_h} \eta \rho v  = \int_{\Omega_h^0} \eta^0 \rho^0 v^0,  \\
\label{eq_energy}
  & \int_{\Omega_h} \eta \rho\, e(\rho,p) + \frac{1}{2} \eta \rho v^2 =
	 \int_{\Omega_h^0} \eta^0 \rho^0 e^0 + \frac{1}{2} \eta^0 \rho^0 (v^0)^2, \\
\label{eq_eta_bounds}
  & 0 \leq \eta_i^{\min} \leq \eta_i \leq \eta_i^{\max}\leq 1,  \\
\label{eq_rho_bounds}
  & \rho_i^{\min} \leq \rho_i \leq \rho_i^{\max},  \\
		\label{eq_p_bounds}
  & p_i^{\min} \leq p_i \leq p_i^{\max},
\end{align}
\end{subequations}
where the vectors of local lower and upper bounds
$(\eta^{\min}, \eta^{\max}, \dots)$ are defined in Section~\ref{sec_bounds}.
After optimization, the specific internal energy $e$ is recovered from the
optimized density and pressure $(\rho,p)$ using the equation of state.
The remap thus returns the same set of state fields
$(\eta,\rho,e,v)$ as provided at the input.

\begin{remark}
The Laghos miniapp \cite{Laghos, Dobrev2012},
which serves as the starting simulation code for our remapping extensions,
evolves specific internal energy as a finite element function.
For the pressure-controlled remap, the approach instead retains the
initial energy of each Lagrangian phase as a Q-function $e_q^0$ and
evolves an FE increment $\delta e_h(t)$, initialized by $\delta e_h(0)=0$.
At every quadrature point $x_q$, the energy used by the Lagrangian operator is
\begin{equation}
  e(x_q,t) = e_q^0(x_q) + \delta e_h(x_q,t).
\end{equation}
The Lagrangian phase requires specific internal energy only at quadrature points.
Therefore, if $e_q^0(x_q)$ is initialized from the corresponding FE energy, this
representation gives exactly the same quadrature-point energy as the standard
FE representation throughout the evolution.
Before each remap step, the reconstructed energy $e(x_q,t)$
is available directly as a Q-function;
after the remap, it is reset to the new $e_q^0$, while $\delta e_h$ is set to zero.
\end{remark}

\begin{remark}
The proposed optimization approach requires an explicit equation of state, in
particular an explicit relation $e=e(\rho,p)$ for \eqref{eq_energy} and
the energy reconstruction after the pressure-controlled remap.
Extensions to tabulated equations of state are left for future work.
\end{remark}


\section{Interpolation and initial guess}
\label{sec_interp}

The direct interpolation in physical space used as the initial step in
the employed remap strategy is performed using the GSLIB library \cite{gslib}.
For brevity, the interpolation procedure is presented only for the indicator field $\eta$.
The density and auxiliary pressure at quadrature point values,
as well as the velocity FE function, are interpolated analogously.

Given source and target meshes $\Omega_h^0$ and $\Omega_h$ over
the same physical domain, each discrete state field on $\Omega_h^0$
is evaluated at any prescribed physical point $x$ in$\Omega_h$.
GSLIB provides the following information for any physical point $x$ in $\Omega_h$:
\begin{enumerate}
	\item The source element $K_x^0 \in \Omega_h^0$ containing $x$ and
	      the corresponding local reference coordinate $\hat{x}_0(x)$,
	      satisfying $T(\hat{x}_0(x)) = x_0(x) = x$,
	      where $T$ denotes the reference-to-physical mapping for $K_x^0$.
  \item The interpolated value at the physical position $x$ on $\Omega_h^0$
        is evaluated as $\eta^0(x) = \eta^0(x_0(x))$.
\end{enumerate}
The above functionality is used for the interpolation of both FE functions
and pointwise Q-functions, as described below.
The interpolation achieves optimal sharpness
but does not conserve the mass of $\eta^0$.


\subsection{Interpolating FE functions}
\label{sec_gslib_fe}
The interpolation of an FE function from $\Omega_h^0$ to $\Omega_h$ evaluates the
source field $\eta^0$ at the coordinates of the support nodes of the target mesh.
The resulting values serve as the FE coefficient vector of $\eta$.
These coefficients naturally inherit the local bounds of the sampled source values.
Furthermore, representing higher-order FE functions
in a positive Bernstein basis ensures that all elemental field values
lie within the range of the element coefficients, thus preserving the bounds.


\subsection{Interpolating Q-functions}
\label{sec_gslib_Q}
In general, quadrature point values are defined locally and do not have a
physical area of support or an interpolation.
Interpolation onto a new mesh therefore requires a spatial reconstruction.
Two reconstruction options are considered.
The first extends the source quadrature point values $q^0$ to a
piecewise constant function $\overline{q^0}$ on a refined mesh
$\overline{\Omega_h^0}$, with one element per quadrature point
assigned the corresponding value of $q^0$.
The target values $q$ are then obtained by
evaluating $\overline{q^0}$ at the quadrature points of $\Omega_h$,
following the procedure described in Section~\ref{sec_gslib_fe}.
This piecewise constant option is appropriate for strong shocks and
discontinuities, where a non-oscillatory reconstruction is more important
than high-order accuracy in smooth regions.
The current implementation in MFEM~\cite{MFEM2020} supports this construction
only for source meshes composed of tensor-product elements, namely
quadrilaterals or hexahedra.

The second option is a high-order Gauss-Legendre reconstruction.
When values are stored at Gauss-Legendre quadrature points, they are
interpreted as the nodal coefficients of an elementwise discontinuous $L^2$
polynomial of the same degree, expressed in a Gauss-Legendre basis.
This defines a discontinuous FE function on the source mesh $\Omega_h^0$,
whose nodes coincide with the source quadrature points.
The target Q-function values are then obtained by interpolating this FE
function at the new-mesh quadrature points,
using the standard FE interpolation procedure of Section~\ref{sec_gslib_fe}.
This option retains the high-order approximation of smooth Q-functions.


\section{Definition of admissible bounds}
\label{sec_bounds}

The optimization formulation enforces local bounds
on the state fields through the constraints in
Eqs.~\eqref{eq_eta_bounds}, \eqref{eq_rho_bounds}, and \eqref{eq_p_bounds}.
These bounds are constructed from interpolated values
at the support points of the target mesh.
For FE functions, they are the points associated with the degrees of freedom.
Let $S_K$ denote the set of support points in a target element $K\in\Omega_h$.
For the indicator field, the bounds are computed as
\begin{equation}
	\label{eq_minmax}
	\eta_K^{\min} = \min_{z\in S_K}\eta^*(z), \qquad
	\eta_K^{\max} = \max_{z\in S_K}\eta^*(z).
\end{equation}
Thus, the indicator is bounded by the range of its direct interpolation
within each new-mesh element.
In the case of $\eta_K^{\max}<\varepsilon$, for a small threshold $\varepsilon$,
the element is treated as empty, and both indicator bounds and
interpolated indicator values are set to zero.

Other local-bound definitions are possible. For example, one could construct
bounds from a wider stencil of elements on the initial mesh. The
new-mesh elementwise construction above is used because it limits support propagation
in a remap to at most one layer of new-mesh elements.

Density and auxiliary pressure are represented as Q-functions
in the present hydrodynamic formulation.
As such, their bounds are constructed analogously to the indicator bounds,
but with an indicator-aware restriction.
For each nonempty element $K$, let $Q_K^+ = {q\in Q_K:\eta^*(q)>\varepsilon}$
denote the set of active quadrature points.
The bounds for $f\in{\rho,p}$ are then defined by
\begin{equation}
	f_K^{\min} = \min_{q\in Q_K^+} f^*(q), \qquad
	f_K^{\max} = \max_{q\in Q_K^+} f^*(q).
\end{equation}
These bounds are assigned to every quadrature point of $K$,
including points where $\eta^*(q)\leq\varepsilon$.
If $K$ is empty, the density and pressure bounds are both set to zero.
Consequently, the density and pressure are supported on the same set of active
new-mesh elements as the indicator; this elementwise synchronization is
required for a consistent coupled remap.
No pointwise bounds are imposed on velocity, which is retained as its
interpolated value during the optimization.

\begin{remark}
State variables with identical lower and upper bounds are fixed and can therefore
be excluded from the set of optimization variables to reduce computational cost.
\end{remark}


\section{Optimization methods}
\label{sec_optimization}

This section summarizes the basic concepts of the two optimization algorithms
used in this work, the IPM implemented in HiOp \cite{Petra2019}
and the Bregman projection algorithm.
Although alternative algorithms, such as the Method of Moving Asymptotes
\cite{Svanberg2002}, could also be employed, IPM is selected for its broad
applicability and established convergence theory, while the entropic Bregman
projection approach is selected for its natural handling of design-variable bounds.


\subsection{Interior point method (IPM)}
\label{sec_opt_hiop}

HiOp’s interior-point implementation is designed for large-scale, nonlinear, nonconvex optimization problems with simple bounds and a modest number of global equality and inequality constraints. It combines a primal-dual filter interior-point framework~\cite{waecther_05_ipopt0} with a quasi-Newton solution strategy~\cite{HiOpRepo}. The Hessian of the Lagrangian is approximated using limited-memory BFGS updates based on first-order derivative information, avoiding the computationally expensive evaluation of exact second derivatives. The step computation is embedded in the filter line-search framework~\cite{waecther_05_ipopt}, so trial points are accepted based on a two-criterion decrease mechanism that balances infeasibility and barrier-objective reduction~\cite{waecther_05_ipopt}, rather than through a penalty function alone. This approach preserves the global convergence features of filter interior-point techniques while reducing the derivative and linear-algebra costs associated with full Newton IPMs~\cite{Petra2019}.

\subsection{Entropic Bregman Projection (BP) algorithm}
\label{sec_opt_bp}

The minimization problem is re-formulated as:
\begin{equation}\label{eq:BP}
	\begin{aligned}
		 & \min_{\eta, \rho, p} ~ \int_{\Omega_h}\Big(D_\varphi(\eta,  \eta^*) + D_\varphi(\rho,  \rho^*) + D_\varphi(p,  p^*)+D_\varphi(v, v^*)\Big) \\
		 & \quad\text{subject to:}\quad \eqref{eq_vol}\text{-}\eqref{eq_energy}.
	\end{aligned}
\end{equation}
where $D_\varphi:(f,g)\mapsto \varphi(f)-\varphi(g)-\nabla \varphi(g)\cdot (f-g)$ denotes the Bregman divergence \cite{Bregman1967} associated with the generalized Fermi-Dirac binary entropy $\varphi$, which encodes the bound constraints of the original problem:
\begin{equation*}
	\varphi(f)
	=\begin{cases}
		\begin{aligned}          & \Big((f-f^{\textrm{min}})\ln (f-f^{\textrm{min}})         \\
                         & \quad + (f^{\textrm{max}}-f)\ln (f^{\textrm{max}}-f)\Big)\end{aligned} & \text{ when }f^{\textrm{min}}\leq  f \leq f^{\textrm{max}}, \\
		+\infty                                                                                                                       & \text{ otherwise.}
	\end{cases}
\end{equation*}
Here, the bounds $f^{\textrm{min}}$ and $f^{\textrm{max}}$ for an optimized variable
$f\in\{\eta,\rho,p\}$ are defined at the quadrature point value.
The inverse derivative of $\varphi$ is a scaled sigmoid function $\sigma(\psi_f)=f^{\min}+(f^{\max}-f^{\min})/(1+\exp(-\psi_f))$.

Since the constraints are non-convex, a local minimizer may not be the global minimizer. However, assuming the interpolated solutions are close to the globally conservative ones, local minima near them suffice for the present purpose.

To find a local minimizer of \eqref{eq:BP}, Dykstra's projection algorithm \cite{Bauschke2000} is employed. It computes the projection onto an intersection $C=\bigcap_{i=1}^m C_i$ of closed convex sets by cyclically projecting onto each $C_i$ while accumulating correction terms that drive the iterates toward the projection onto $C$ rather than onto any single $C_i$.
This algorithm is guaranteed to be convergent when each $C_i$ is closed and convex \cite{Bauschke2000}.
One set $C_i$ is associated with each global constraint in \eqref{eq_problem}.
The volume constraint \eqref{eq_vol} defines an affine (hence convex) set, whereas the mass, momentum, and energy constraints \eqref{eq_mass}-\eqref{eq_energy} are nonlinear in the primal variables and therefore non-convex.
Each non-convex $C_i$ is replaced at each projection step by its tangent affine subspace at the current iterate, and a Bregman-divergence Dykstra projection algorithm is applied to the resulting collection of affine sets.

\begin{algorithm}[ht]
\caption{Dykstra projection}\label{alg:dykstra}
\small
\begin{algorithmic}
\Require $u^*:=(\eta^*, \rho^*, p^*)$, $v^*$, $\texttt{tol}$, $n_{\texttt{lin}}$
\Ensure $u=:(\eta, \rho, p)$ satisfying \eqref{eq_vol}-\eqref{eq_energy} up to tolerance $\texttt{tol}$, and $e$.
\State Let $r_1,\dots,r_m$, $m=3+d$, be the residual functions of \eqref{eq_vol}, \eqref{eq_mass}, \eqref{eq_energy}, and \eqref{eq_momentum}.
\State $\psi \gets \nabla \varphi(u^*)$, with clipping if necessary
\State Initialize $q_1,\dots,q_m$ with zero vectors.
\While{$\max_{i} |r_i(\sigma(\psi))|>\texttt{tol}$}
    \For{$i=1,\dots,m$}
        \State $y \gets \psi + q_i$, \quad $\psi \gets y$
        \For{$j=1,\dots,n_{\texttt{lin}}$} \Comment{tangent projections onto $C_i$}
            \State $u \gets \sigma(\psi)$, \quad $g \gets \nabla r_i(u)$
            \State Find $\mu\in\mathbb{R}$ such that
            \begin{equation}
            \label{eq:proj}
                \int_{\Omega_h} g\cdot \sigma(\psi + \mu g) = \int_{\Omega_h} g\cdot u - r_i(u)
            \end{equation}
            \State $\psi \gets \psi + \mu g$
            \If{$|r_i(\sigma(\psi))| \leq \texttt{tol}$} \textbf{break} \EndIf
        \EndFor
        \State $q_i \gets y - \psi$
    \EndFor
\EndWhile
\State $u\gets \sigma(\psi)$
\State $e \gets e(\rho, p)$ pointwise at the quadrature points
\State \Return $u$, $e$
\end{algorithmic}
\end{algorithm}

Note that $\varphi$ is singular at boundaries $f^{\text{min}}$ and $f^{\text{max}}$, which lead to numerical instabilities. Instead, motivated by \cite{keith2026proximal}, the latent variable $\psi_f=\sigma^{-1}(f)$ is employed.
Then the Bregman projection onto a single linearized equality constraint reduces to the pointwise update
\begin{equation*}
	\psi_f \leftarrow \psi_f + \mu g_f, \quad f \leftarrow \sigma(\psi_f),
\end{equation*}
where $g_f$ is the gradient of the global nonlinear constraints and $\mu \in \mathbb{R}$ is the Lagrange multiplier determined by enforcing the constraint itself.
This can be done by solving a strictly monotone scalar equation \eqref{eq:proj} described in \cite{kim2025simple,keith2025remap}, and can be efficiently solved using the Illinois algorithm \cite{Dowell1971}.
The bound constraints \eqref{eq_eta_bounds}-\eqref{eq_p_bounds} are inherited automatically through the effective domain of $\varphi$, so no separate projection is required for them.
The Dykstra sweep (outer iteration) is then repeated until the original nonlinear constraints \eqref{eq_vol}-\eqref{eq_energy} are satisfied within a prescribed tolerance. The algorithm is described in Algorithm~\ref{alg:dykstra}.
Since this algorithm is a fixed-point algorithm, Anderson acceleration further accelerates it \cite{anderson}.


\begin{remark}
    The objective function of \eqref{eq:BP} is proper only when the initial interpolant $(\eta^*, \rho^*, p^*)$ satisfies the bound constraints \eqref{eq_eta_bounds}-\eqref{eq_p_bounds} strictly.
    If it does not, a clipped interpolant is used.
	In practice, this involves clipping the latent interpolant $\psi^* = \sigma^{-1}(\eta^*)$ before solving the correction equation.
\end{remark}

\begin{remark}
	Within each Dykstra sweep, the tangent-space projection onto a non-convex $C_i$ can be iterated with a predefined maximum number of linearization steps, $n_{\texttt{lin}}$. The tangent affine subspace is reinitialized to better approximate the projection onto the original non-convex set. The results use $n_{\texttt{lin}}=6$.
\end{remark}


\section{Numerical results}
\label{sec_results}

This section demonstrates the accuracy, property preservation, and
robustness of the proposed remap.
Smooth scalar and hydrodynamics tests first verify
the expected high-order convergence rates.
The Three-Shape benchmark then tests coupled hydrodynamic fields with sharp
interfaces, demonstrating the preservation of bounds and global conservation,
correct interface behavior, and consistent support among the remapped fields.
Finally, the standard triple-point hydrodynamics benchmark
in two and three dimensions demonstrates robustness for
full ALE calculations with shocks and vortical flow.

All results are obtained through the MFEM finite element library
\cite{MFEM2020, MFEM2024}, combining its mesh optimization
miniapps \cite{TMOP2019}, Remhos \cite{Remhos}, and Laghos \cite{Laghos}.


\subsection{Convergence on a smooth scalar problem}
\label{sec_result_scalar}

\begin{figure}[pos=htbp]
\centering
\includegraphics[height=0.25\textwidth]{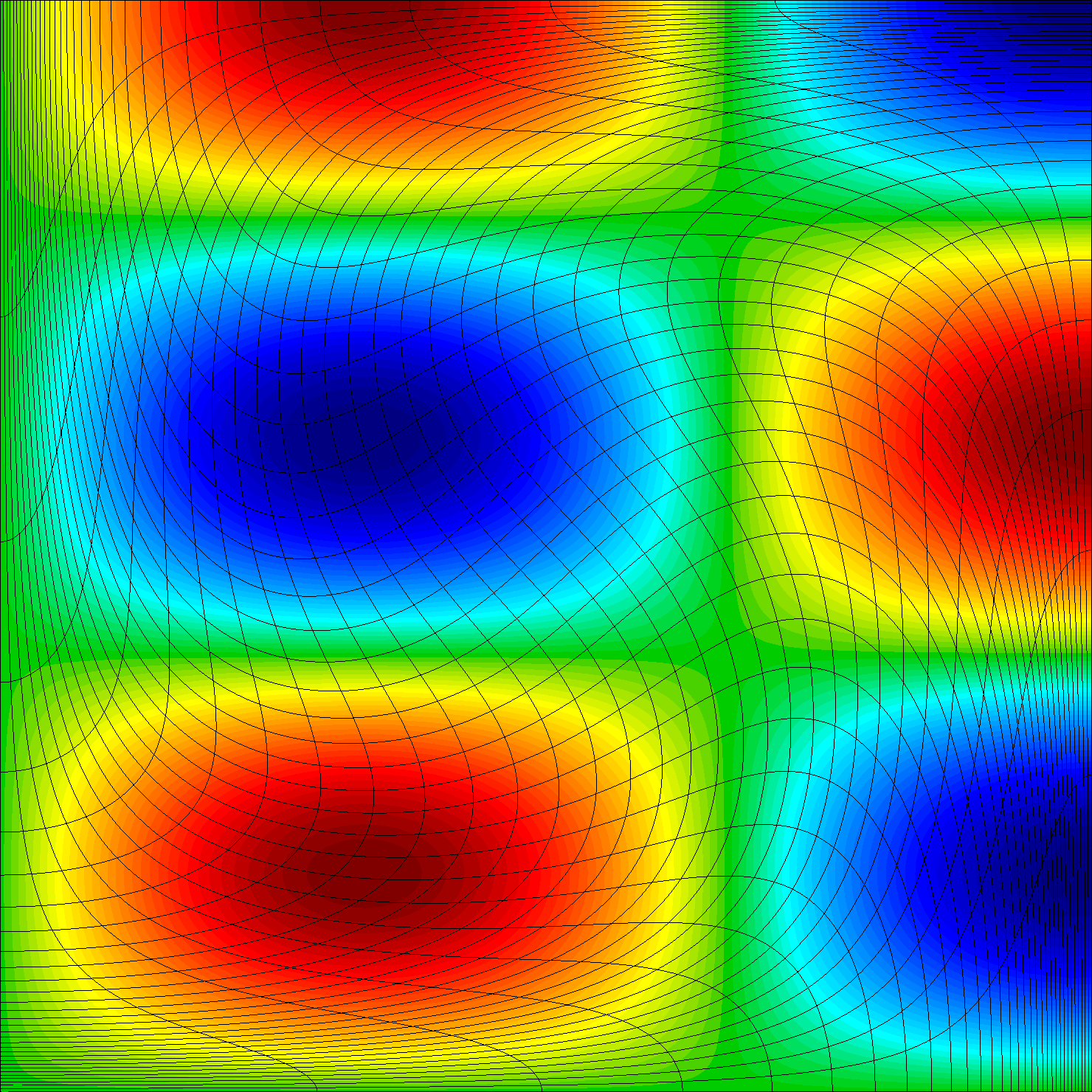}\hfill
\includegraphics[height=0.25\textwidth]{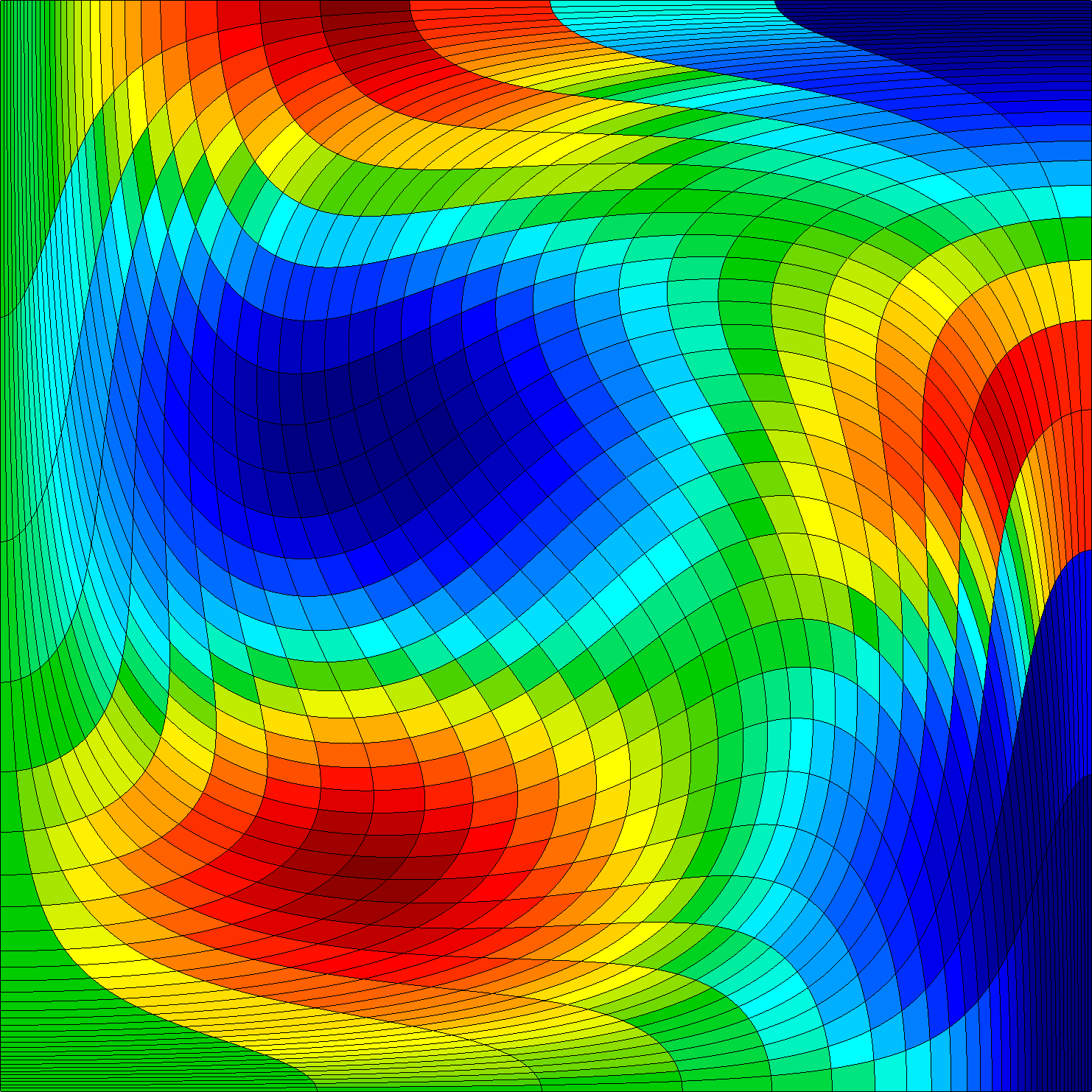}\hfill
\includegraphics[height=0.25\textwidth]{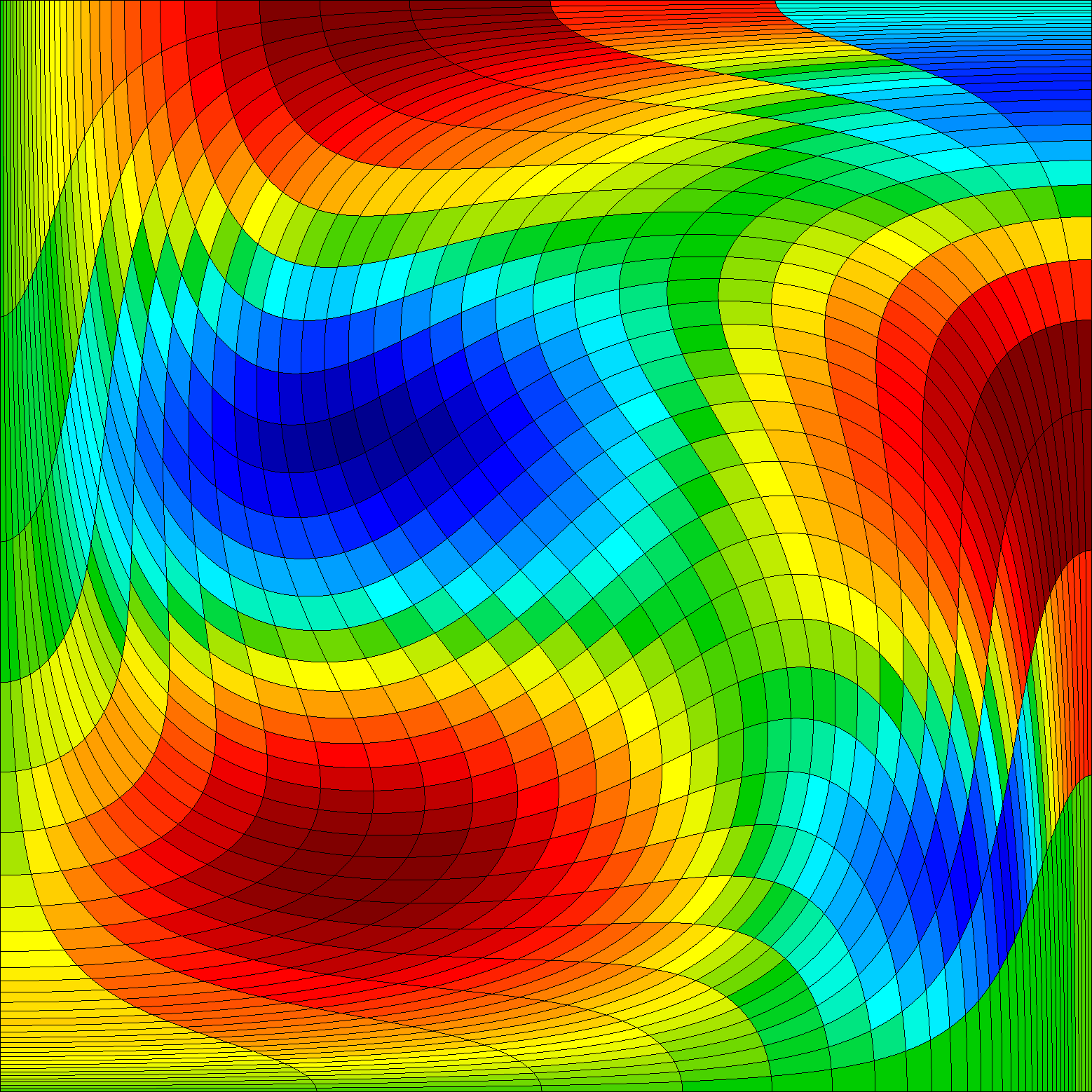}\hfill
\includegraphics[height=0.25\textwidth]{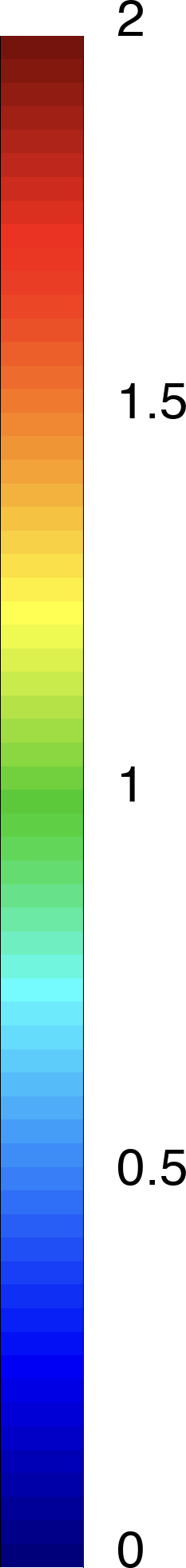}
\caption{FE function result with three uniform refinements.
	       From left to right: final optimized solution on the deformed mesh,
	       minimum bound function, and maximum bound function.
	       The corresponding Q-function plots are visually similar.}
\label{fig:scalar_smooth}
\end{figure}

This test is used to verify the accuracy of the interpolation-optimization
remap independently of the hydrodynamics update.
The smooth scalar function
\[
	u(x,y) = 1 + \sin(3\pi x/2)\sin(5\pi y/2)
\]
is defined on the unit square.
This source function is interpolated directly at the target support points
on a deformed mesh and is then corrected by the optimization.
Polynomial degrees $Q_2$ and $Q_3$ are considered with one to
six uniform refinements, corresponding to mesh sizes $h=1/8,\ldots,1/256$.
For the FE remap, the target support points are the degrees of freedom of a
positive-Bernstein representation.
For the Q-function remap, they are the Gauss-Legendre quadrature points;
the resulting values are interpreted as the coefficients of an elementwise
Gauss-Legendre $L^2$ function of the same degree for the GSLIB interpolation.
Figure~\ref{fig:scalar_smooth} shows a representative FE function result
for the $Q_3$ run with three uniform refinements.

Tables~\ref{tab_remap_fe} and~\ref{tab_remap_qf} report the $L^1$ error
$E_1$, the rate $\log_2(E_{h}/E_{h/2})$, and the absolute mass defects before
and after the optimization, denoted by subscripts $I$ and $O$.
The target mass is $M_0=1+1/(3.75\pi^2)\approx 1.027019$.
The initial interpolation already gives a high-order approximation,
as shown by the convergence rates, but it does not
in general preserve the target mass exactly.
The optimization step changes the interpolated values to restore this mass:
the $|\Delta M|_I$ columns show the mass defect after interpolation,
while the $|\Delta M|_O$ columns show the defect after optimization.
The main point of the test is that this correction preserves
the expected accuracy, with third-order convergence for $Q_2$ and fourth-order
convergence for $Q_3$, while reducing the mass error to the prescribed
$10^{-13}$ tolerance or near roundoff.
All reported tables and numerical values in this subsection use the
Bregman-projection method described in Section~\ref{sec_opt_bp}.
The HiOp optimization method produces similar convergence behavior.



\begin{table}[H]
	\centering
	\caption{Convergence of the FE function remap.}
	\label{tab_remap_fe}
	\footnotesize
	\setlength{\tabcolsep}{2.5pt}
	\begin{tabular}{c|c|cccc|cccc}
		\hline
		  &     & \multicolumn{4}{c|}{$Q_2$}
		        & \multicolumn{4}{c}{$Q_3$} \\
		$r$ & $h$ & $E_1$ & rate & $|\Delta M|_I$ & $|\Delta M|_O$
		    & $E_1$ & rate & $|\Delta M|_I$ & $|\Delta M|_O$ \\
		\hline
		1 & 1/8   & 2.30E-2 & --   & 3.18E-4  & 6.79E-14 & 5.61E-3 & --   & 4.19E-4  & 4.00E-15 \\
		2 & 1/16  & 4.16E-3 & 2.47 & 4.92E-6  & 2.89E-15 & 6.39E-4 & 3.13 & 1.30E-5  & 4.44E-16 \\
		3 & 1/32  & 5.75E-4 & 2.86 & 2.13E-7  & 3.11E-15 & 8.33E-5 & 2.94 & 1.28E-7  & 2.89E-15 \\
		4 & 1/64  & 9.16E-5 & 2.65 & 4.45E-9  & 3.55E-15 & 5.18E-6 & 4.01 & 3.41E-9  & 6.22E-15 \\
		5 & 1/128 & 1.17E-5 & 2.97 & 7.71E-11 & 1.91E-14 & 3.30E-7 & 3.97 & 6.07E-11 & 3.42E-14 \\
		6 & 1/256 & 1.48E-6 & 2.98 & 1.25E-12 & 9.33E-15 & 2.10E-8 & 3.97 & 1.01E-12 & 1.55E-14 \\
		\hline
	\end{tabular}
\end{table}


\begin{table}[H]
	\centering
	\caption{Convergence of the Q-function remap.}
	\label{tab_remap_qf}
	\footnotesize
	\setlength{\tabcolsep}{2.5pt}
	\begin{tabular}{c|c|cccc|cccc}
		\hline
		  &     & \multicolumn{4}{c|}{$Q_2$}
		        & \multicolumn{4}{c}{$Q_3$} \\
		$r$ & $h$ & $E_1$ & rate & $|\Delta M|_I$ & $|\Delta M|_O$
		    & $E_1$ & rate & $|\Delta M|_I$ & $|\Delta M|_O$ \\
		\hline
		1 & 1/8   & 1.17E-3 & --   & 2.93E-4  & 8.88E-16
		          & 6.77E-5 & --   & 6.32E-6  & 8.88E-16 \\
		2 & 1/16  & 1.41E-4 & 3.05 & 7.46E-6  & 4.89E-15
		          & 4.09E-6 & 4.05 & 3.56E-7  & 2.22E-16 \\
		3 & 1/32  & 1.73E-5 & 3.03 & 5.99E-7  & 2.66E-15
		          & 2.53E-7 & 4.02 & 4.12E-9  & 2.00E-15 \\
		4 & 1/64  & 2.20E-6 & 2.98 & 3.33E-9  & 6.00E-15
		          & 1.57E-8 & 4.00 & 7.27E-11 & 5.33E-15 \\
		5 & 1/128 & 2.74E-7 & 3.01 & 5.73E-10 & 2.22E-14
		          & 9.83E-10 & 4.00 & 5.93E-13 & 2.58E-14 \\
		6 & 1/256 & 3.42E-8 & 3.00 & 7.78E-11 & 1.33E-15
		          & 6.15E-11 & 4.00 & 1.11E-13 & 2.82E-14 \\
		\hline
	\end{tabular}
\end{table}


\subsection{Two dimensional Taylor-Green vortex}
\label{sec_result_tg}

The two-dimensional Taylor-Green vortex problem is used to assess whether the
conservative remap correction preserves
high-order convergence for a smooth hydrodynamics solution.
The underlying Lagrangian phase is high-order accurate \cite{Dobrev2012}.
In this test, the Lagrangian evolution is interrupted twice for mesh
optimization and remapping.
The remap must therefore be high-order accurate and couple consistently with
both the mesh optimization and Lagrangian phases.
High-order interpolation alone recovers the expected approximation rates.
The main question is therefore whether the subsequent
corrections that enforce the conservation of mass, momentum, and total energy
remain sufficiently small to maintain those rates over mesh refinement.

The initial conditions and Lagrangian evolution
for this problem are described in \cite[Section 8.1]{Dobrev2012}.
The simulation advances to final time $T = 1.25$, with
two ALE steps applied at times $t = 0.5$ and $t = 1.0$.
As specified by this test setup, errors at $T=1.25$ are
measured relative to the initial state.
During each ALE step, mesh optimization restores the
initial uniform Cartesian configuration at $t = 0$.
If both the Lagrangian evolution and the remap are exact, the hydrodynamic
state immediately after each ALE step would coincide with the initial condition,
effectively restarting the problem.

The simulations use the high-order Gauss-Legendre reconstruction from
Section~\ref{sec_gslib_Q} for the density and auxiliary pressure Q-functions.
This choice is particularly important for the pressure, as the
piecewise constant reconstruction would limit the observed convergence rate.
All simulations in this subsection use the Bregman projection method described
in Section~\ref{sec_opt_bp}; the HiOp optimization method produces similar
convergence behavior.

\begin{figure}[pos=htbp]
	\centering
	\begin{minipage}[c]{0.35\textwidth}
		\centering
		\includegraphics[width=\linewidth]{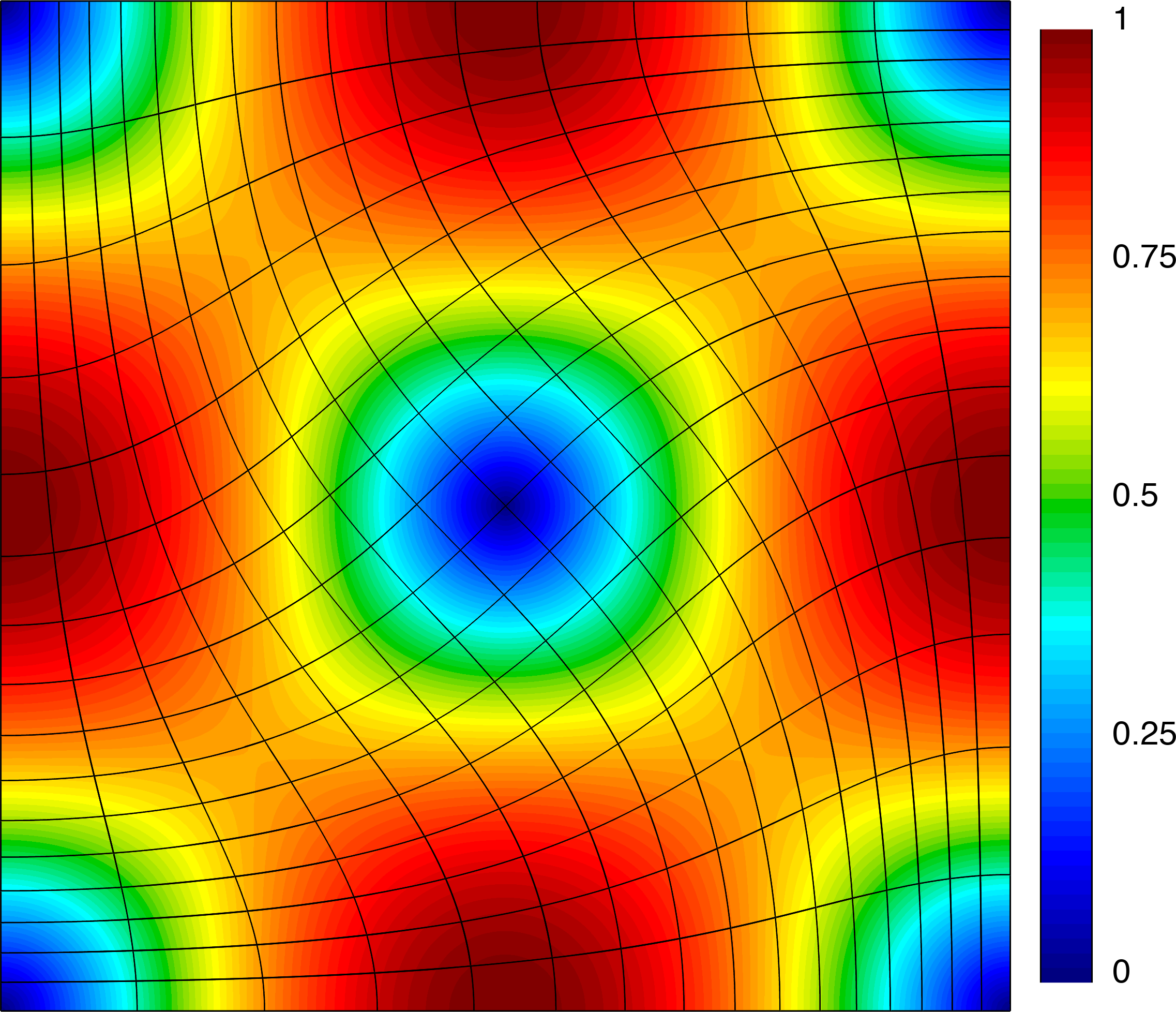}
	\end{minipage}\hfill
	\begin{minipage}[c]{0.55\textwidth}
		\centering
		\includegraphics[width=\linewidth]{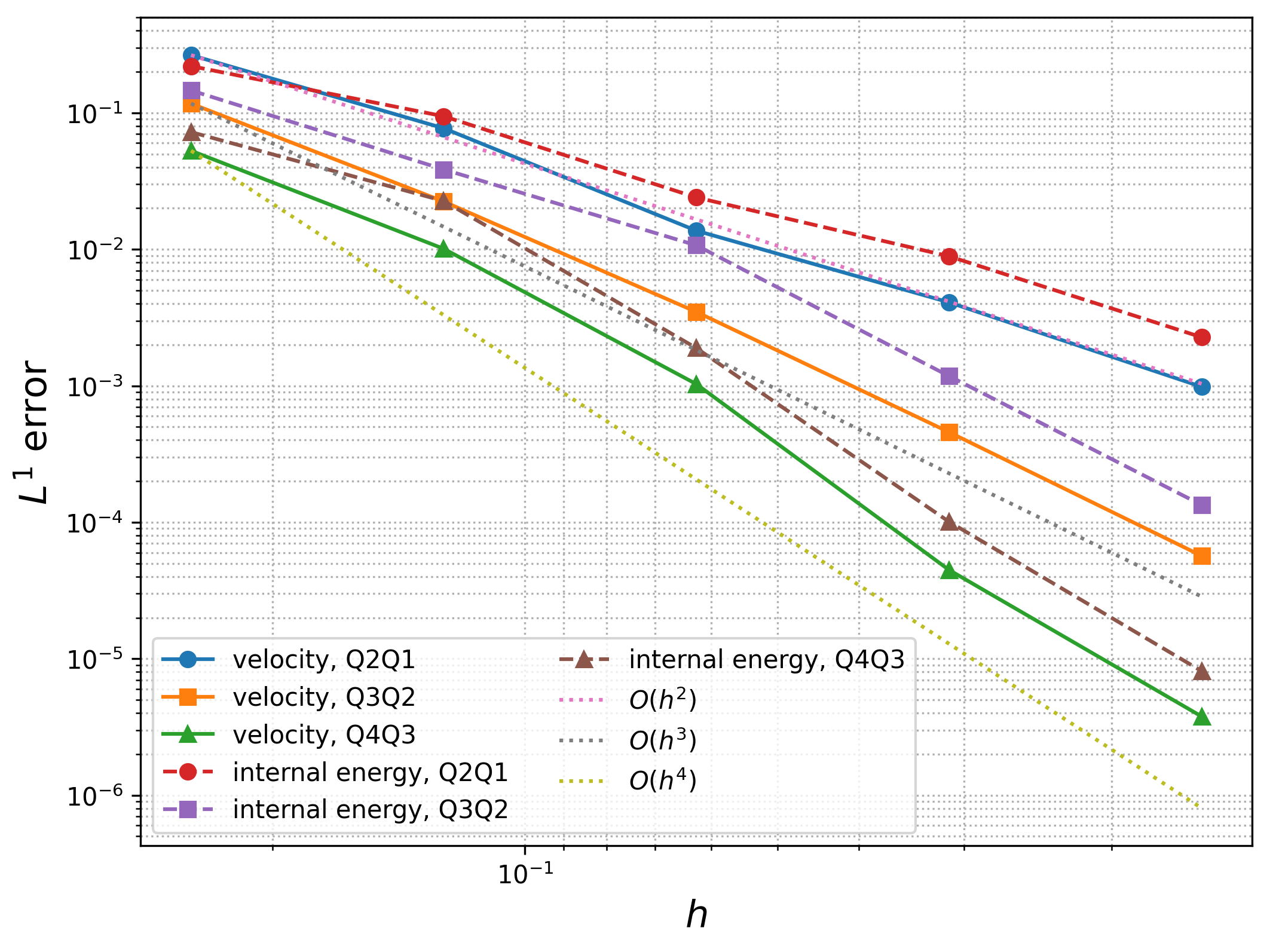}
	\end{minipage}
	\caption{Final $Q_3$ mesh and velocity at $T=1.25$ (left) and
	         convergence rates for different discretizations (right)
	         on the 2D Taylor-Green problem using optimized remap.}
	\label{fig_tg}
\end{figure}

Three discretization pairs, $Q_2Q_1$, $Q_3Q_2$, and $Q_4Q_3$, are tested, with
expected convergence orders of two, three, and four, respectively.
Figure~\ref{fig_tg} confirms these rates.
Convergence beyond second order is uncommon for
schemes combining Lagrangian evolution with ALE remapping.
In addition, the conservation of the total energy is compared
for three setups using the $Q_3Q_2$ discretization in Table~\ref{tab_tg}:
a purely Lagrangian evolution to $t=0.5$,
an interpolation-only remap at $t=0.5$,
and the full, conservative remap including the optimization step at $t=0.5$.
For the optimized remap, the conservation target is
the total energy of the Lagrangian state immediately before remapping.
The optimization is therefore expected to preserve
the existing Lagrangian energy error relative to the exact stationary solution,
up to the optimization tolerance and roundoff.
Table~\ref{tab_tg} confirms this expectation, whereas
interpolation-only remapping degrades the conservation of total energy.
Together with the convergence results, these findings demonstrate that the
optimization step restores conservation while preserving high-order convergence.

\begin{table}[h!]
	\centering
	\caption{Total energy errors at $T=0.5$ for the Q3Q2 Taylor-Green test.}
	\label{tab_tg}
	\begin{tabular}{c|ccc}
		\hline
		$h$  & Lagrangian & Optimized & Interpolated \\
		\hline
		1/4  & 6.61E-7   & 6.60E-7  & 1.58E-3     \\
		1/8  & 5.41E-9   & 5.41E-9  & 4.49E-5     \\
		1/16 & 3.82E-12  & 4.07E-12 & 1.35E-6     \\
		1/32 & 4.83E-14  & 2.51E-14 & 5.98E-7     \\
		1/64 & 2.78E-17  & 4.34E-12 & 1.45E-8     \\
		\hline
	\end{tabular}
\end{table}


\subsection{Two dimensional Three-Shape benchmark}
\label{sec_result_shapes}

The two-dimensional Three-Shape benchmark
introduced in \cite{LeVeque1996} is a standard test of
remap accuracy for complex geometries and discontinuous indicator fields.
This example extends the usual indicator-function test to include density,
pressure, specific internal energy, and velocity fields, assessing the
remapping of the coupled hydrodynamic system rather than a single scalar field.
The pressure-controlled remap updates $(\eta,\rho,p,v)$ and
reconstructs $e$ from the equation of state.

\begin{figure}[pos=htbp]
\centering
  \includegraphics[height=0.25\textwidth]{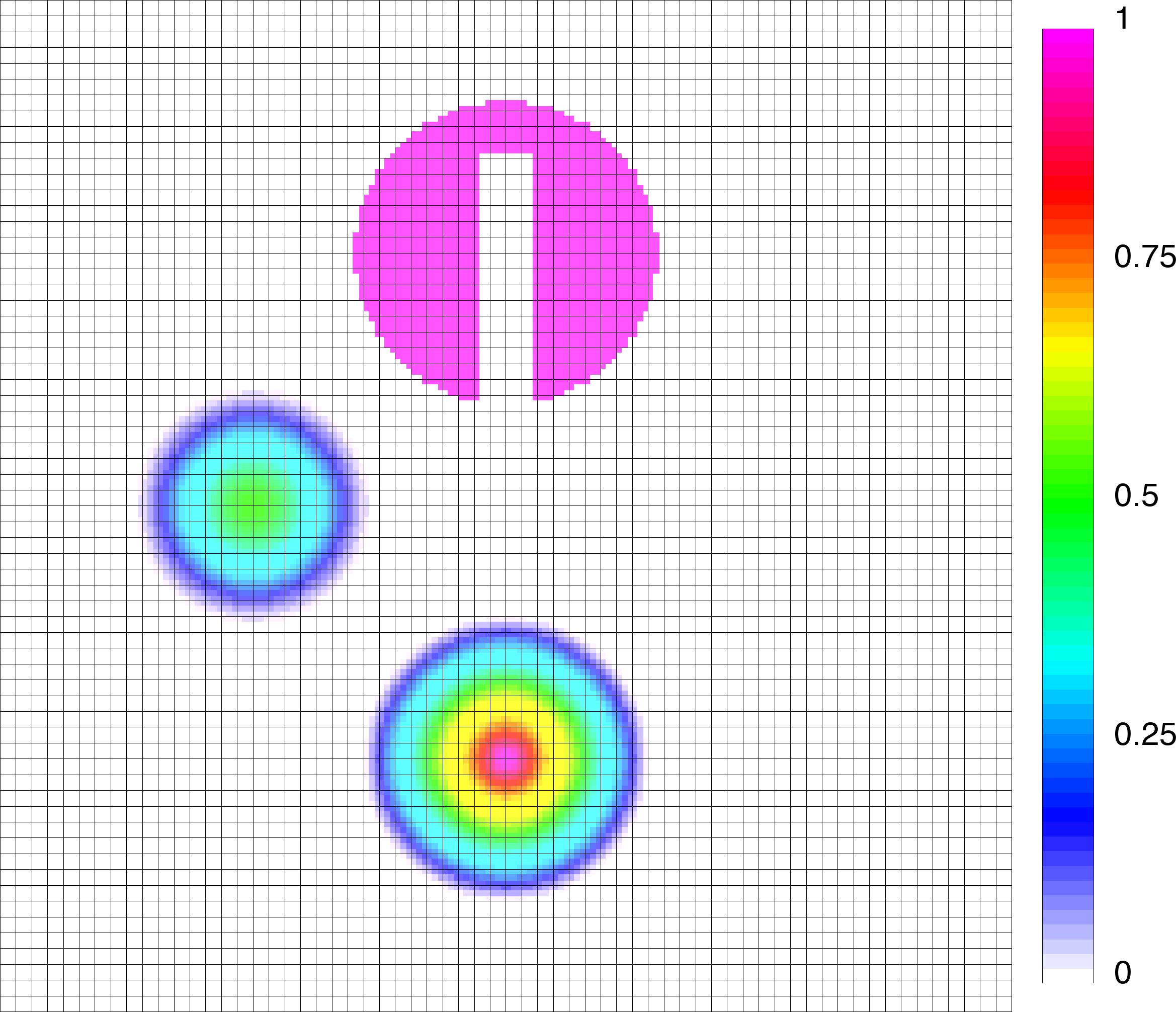} \hfil
  \includegraphics[height=0.25\textwidth]{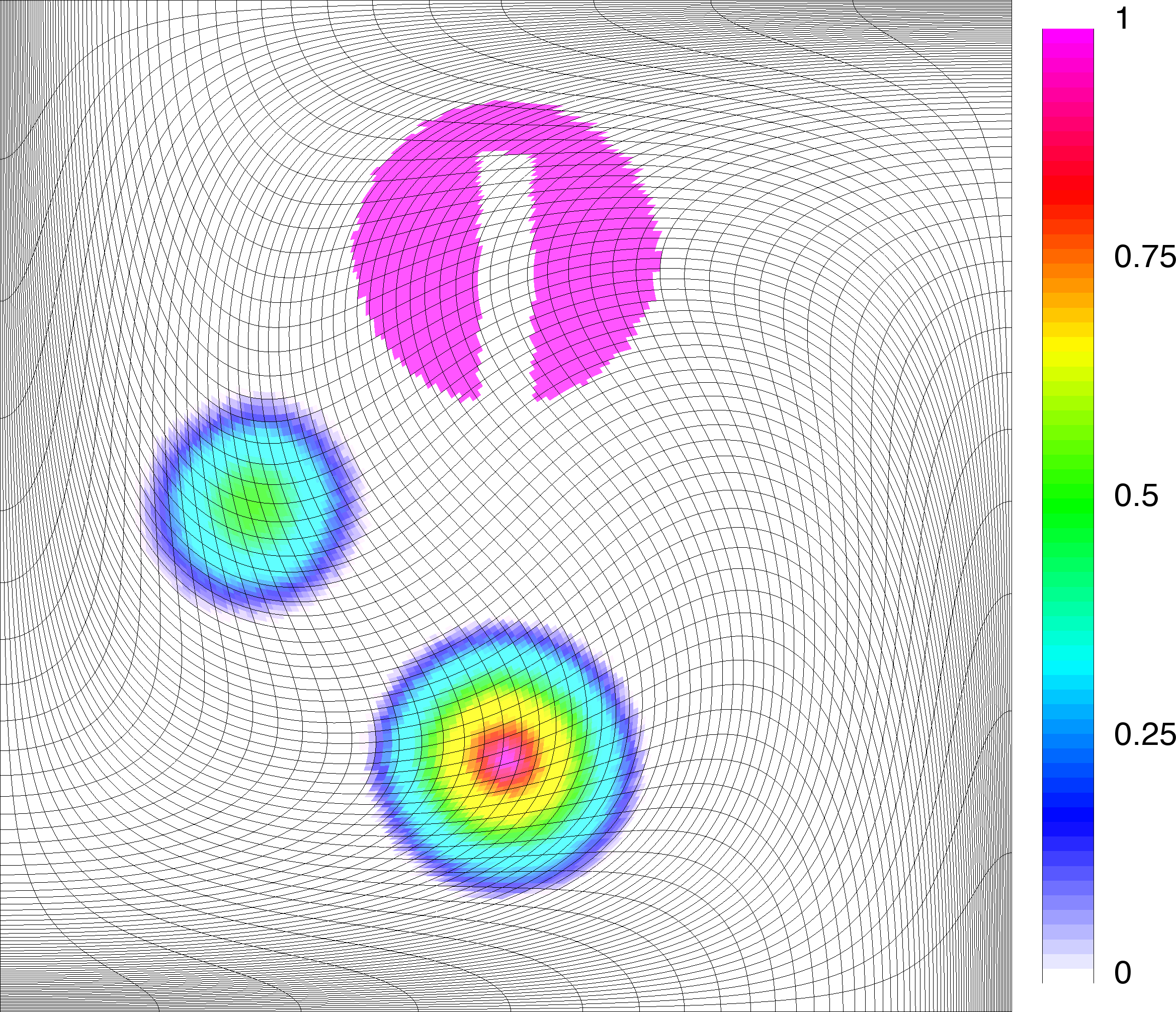} \hfil
  \includegraphics[height=0.25\textwidth]{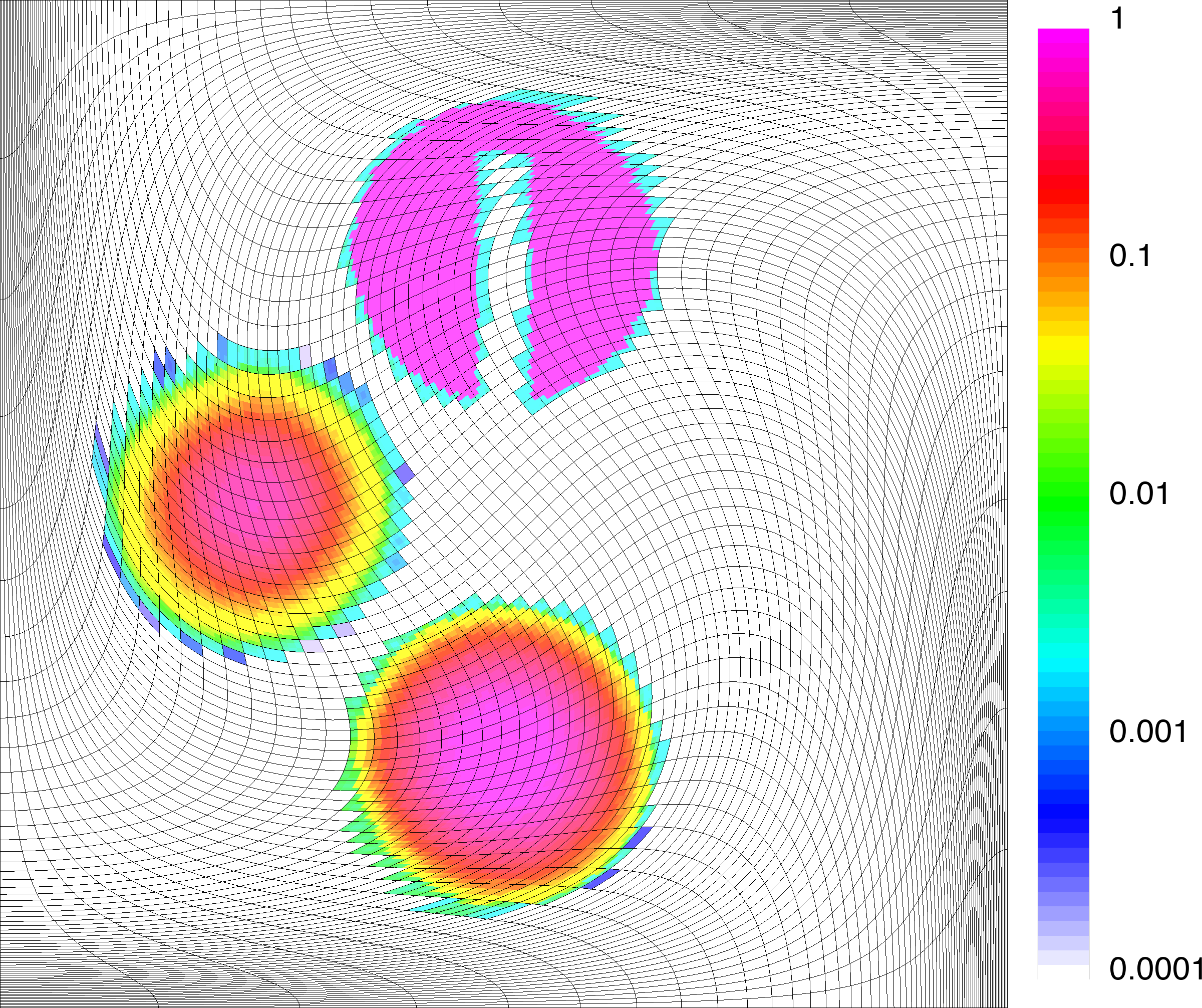}
  \caption{Initial mesh and indicator Q-function (left),
           final mesh and optimized indicator Q-function (middle),
           and log-scale of the optimized indicator Q-function result (right)
           for the Three-Shape benchmark.}
  \label{fig_3s_ind}
\end{figure}

The initial problem is discretized on a
uniform $64 \times 64$ mesh $\Omega_h^0$ with quadratic elements.
The target mesh $\Omega_h$ is generated through a perturbation that evolves
the initially uniform mesh coordinates $x$ over $\tau \in [0,1]$ according to
\[
	\frac{dx}{d \tau} = (\sin(\pi x_0) \cos(\pi x_1), -\cos(\pi x_0) \sin(\pi x_1)),
\]
as illustrated in Figure~\ref{fig_3s_ind}.
The constrained correction for this test is computed with the
HiOp interior-point method described in Section~\ref{sec_opt_hiop};
The BP method yields very similar results.

\begin{figure}[pos=htbp]
\centering
  \includegraphics[height=0.25\textwidth]{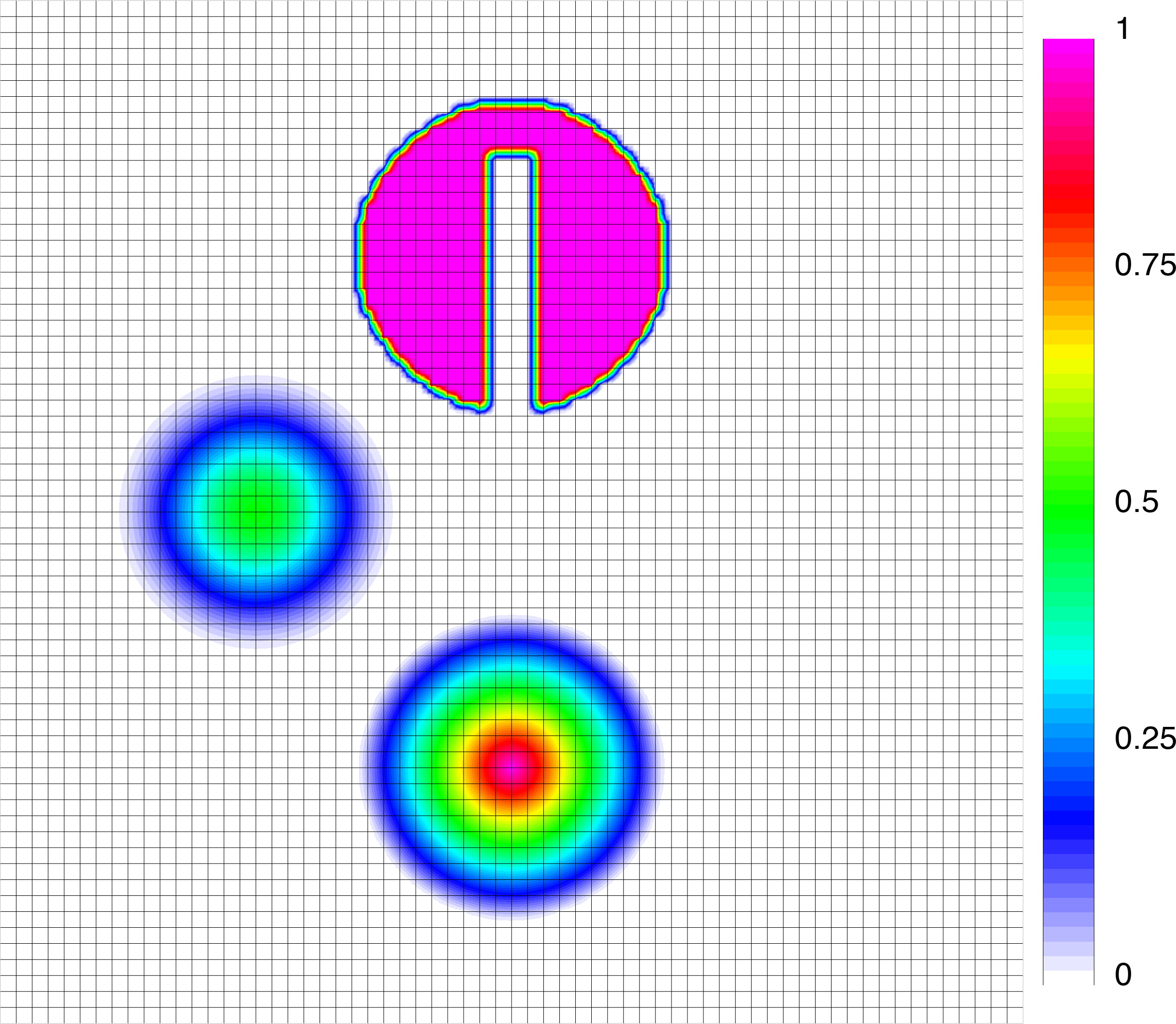} \hfil
  \includegraphics[height=0.25\textwidth]{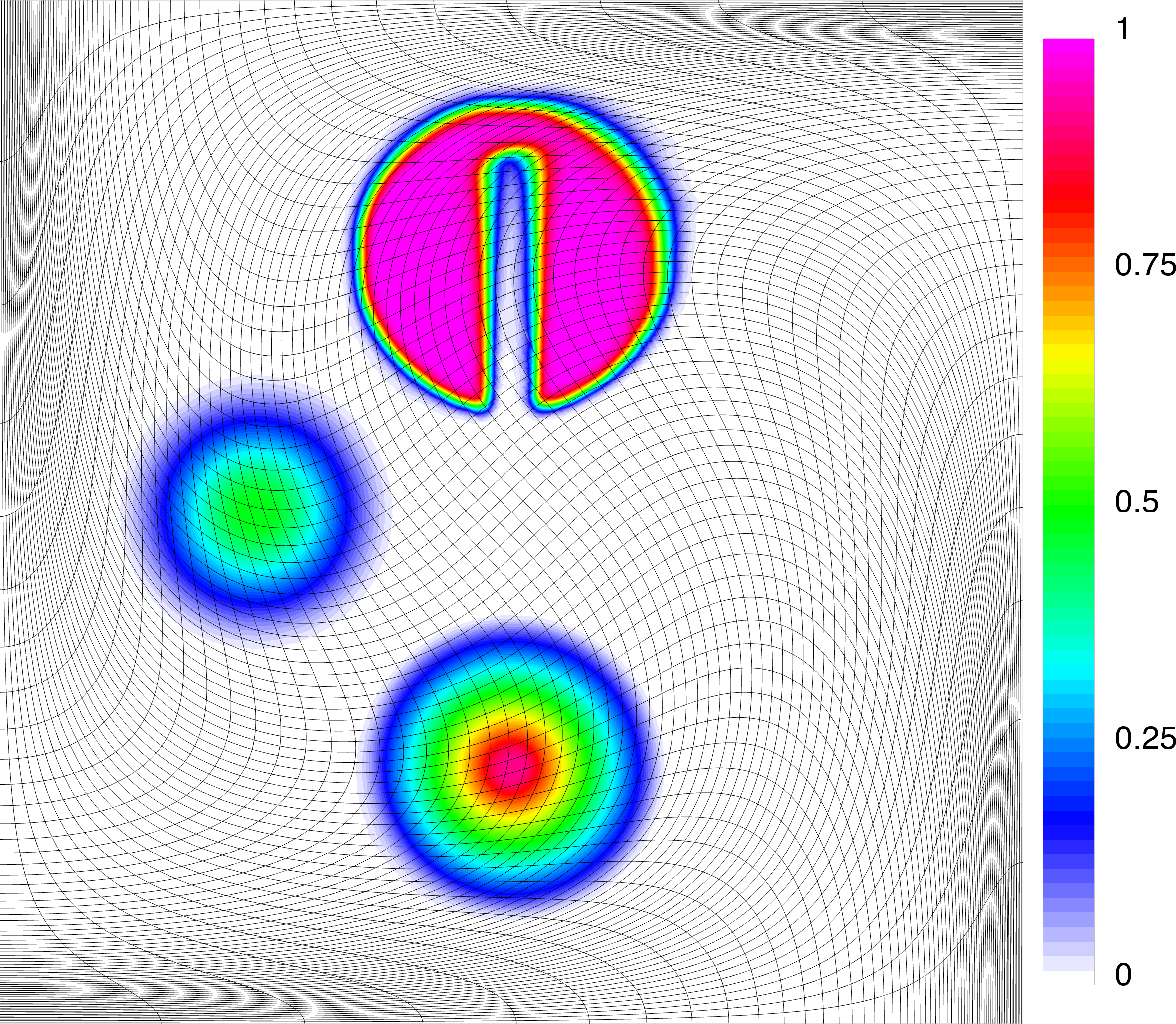} \hfil
  \includegraphics[height=0.25\textwidth]{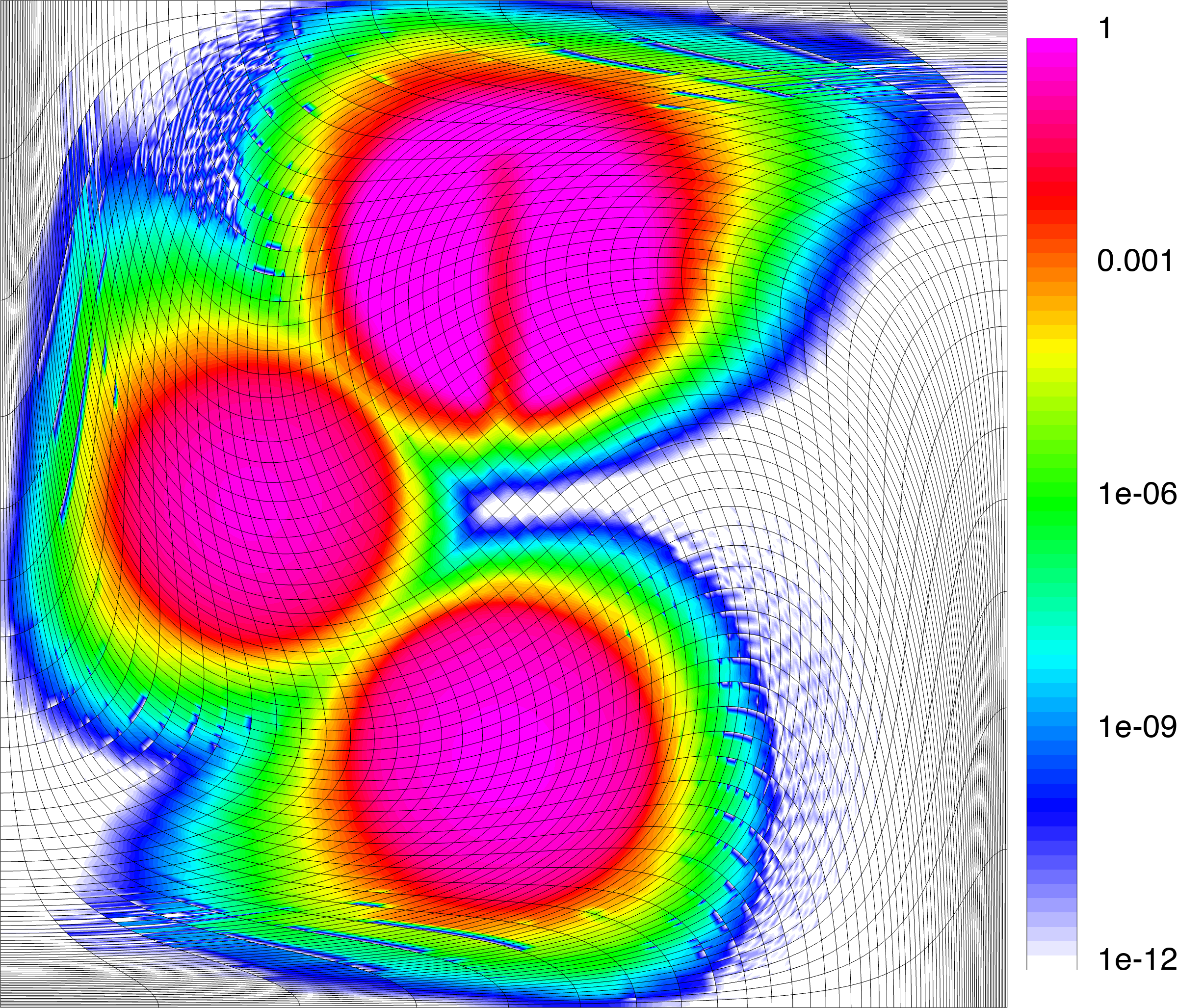}
  \caption{Initial mesh and FE indicator function (left),
           final mesh and advection-remapped FE indicator function (middle),
           and log-scale of the advection-remapped indicator function (right)
           for the Three-Shape benchmark.}
  \label{fig_3s_ind_adv}
\end{figure}

To isolate the diffusion properties of the material indicator,
the indicator field is first considered on its own.
The sharp interface between the material region and its complement tests the
ability to preserve discontinuities with minimal numerical diffusion.
The optimized Q-function remap shown in Figure~\ref{fig_3s_ind} limits the
indicator transition region to approximately one to two elements in width while
maintaining conservation and satisfying the prescribed bounds.
In contrast, the representative advection-based remap shown in
Figure~\ref{fig_3s_ind_adv} visibly diffuses the indicator across the interface.

\begin{figure}[pos=htbp]
\centering
\begin{subfigure}[t]{0.3\textwidth}
  \centering
  \includegraphics[width=\textwidth]{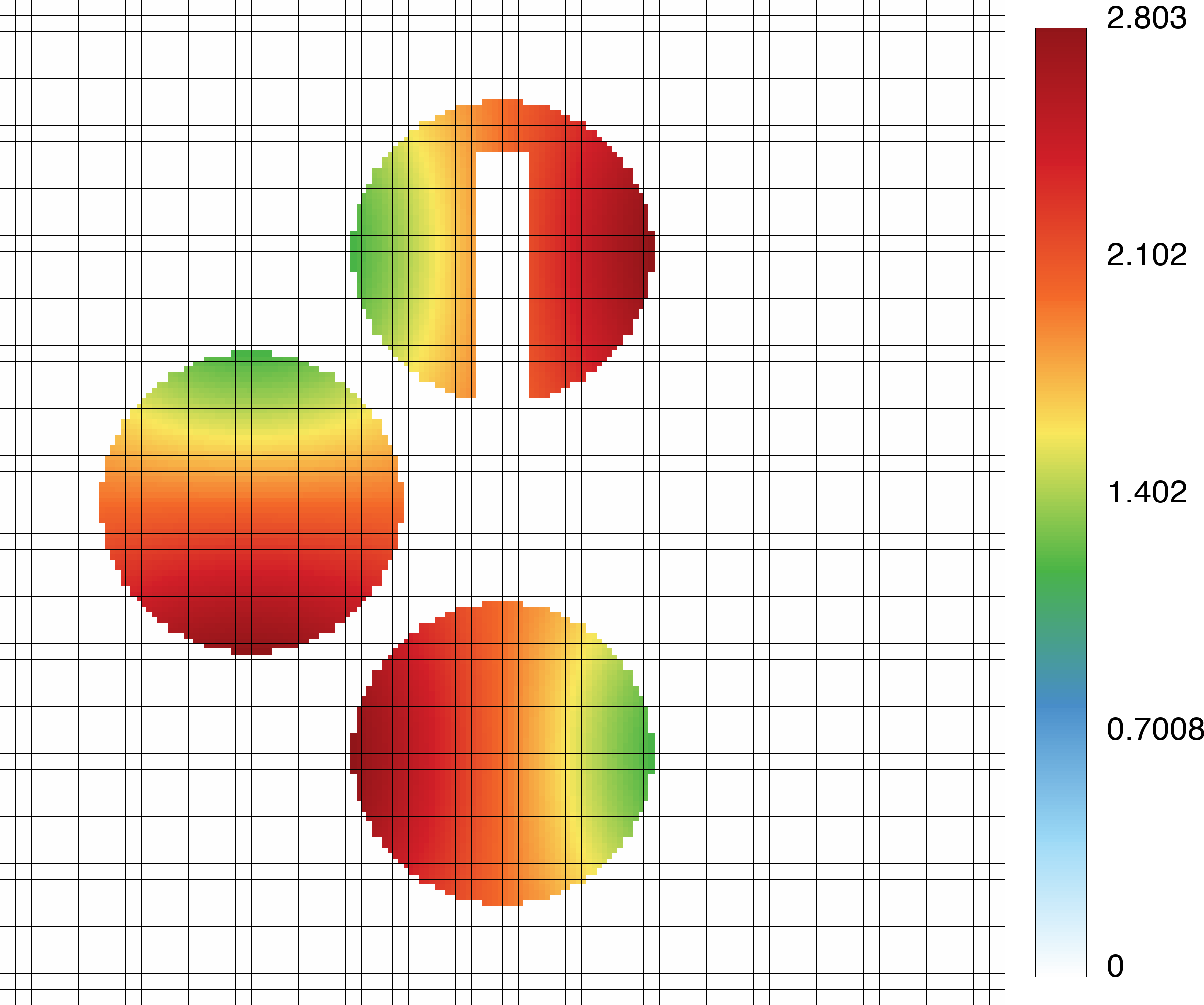}
  \caption{$\rho$, initial}
\end{subfigure}\hfil
\begin{subfigure}[t]{0.3\textwidth}
  \centering
  \includegraphics[width=\textwidth]{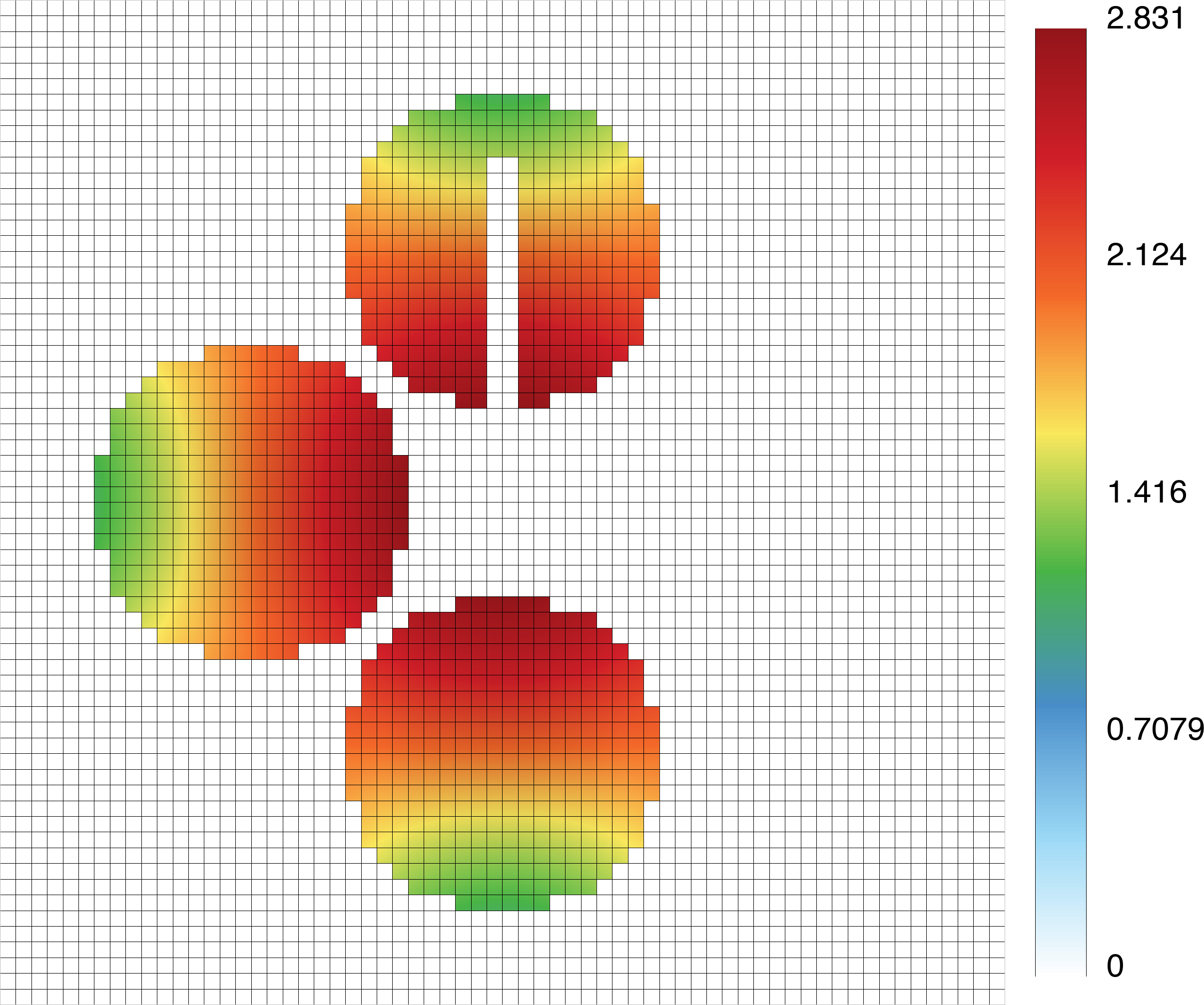}
  \caption{$e$, initial}
\end{subfigure}\hfil
  \begin{subfigure}[t]{0.3\textwidth}
  \centering
  \includegraphics[width=\textwidth]{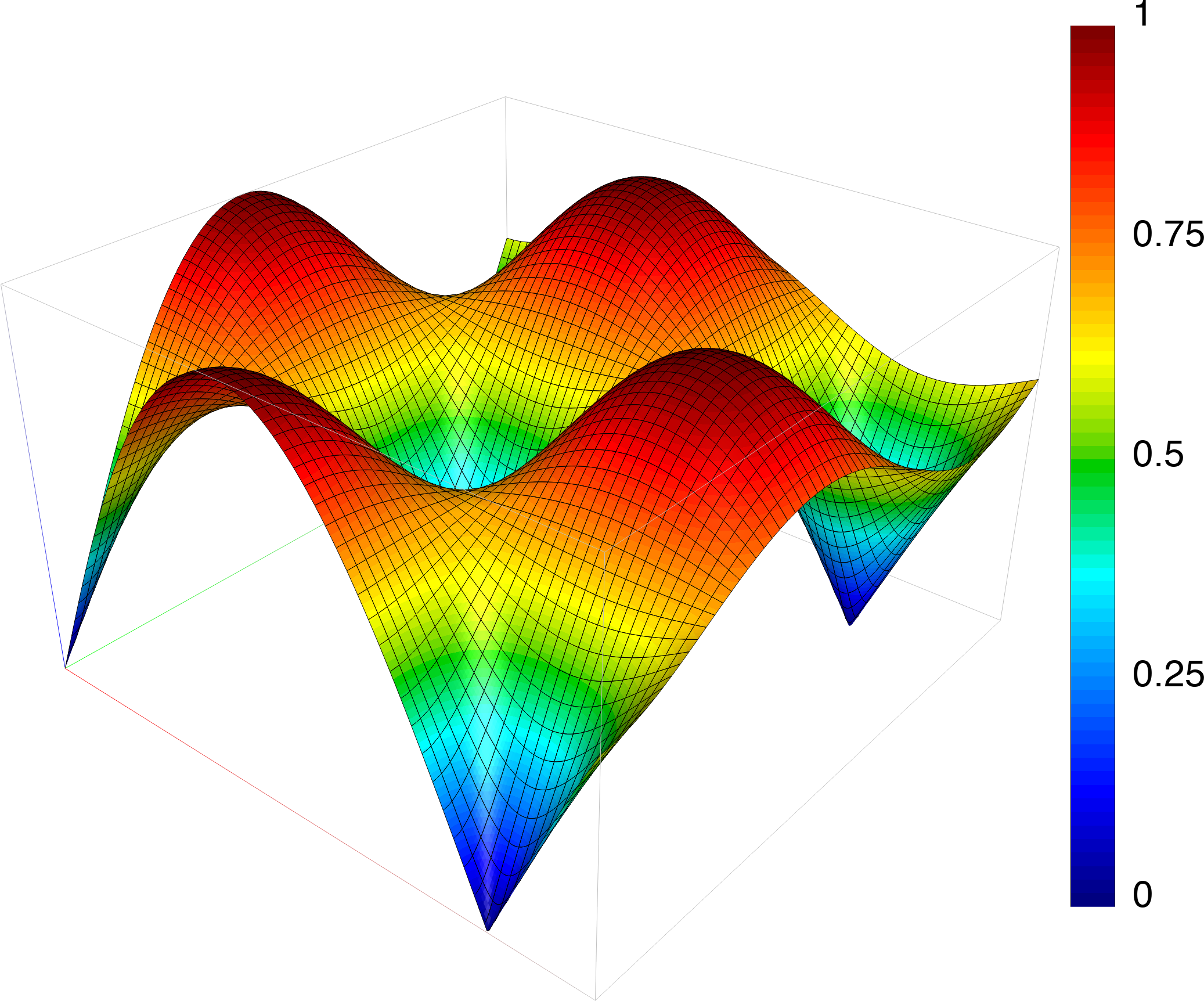}
  \caption{$v$, initial}
\end{subfigure}
\vspace{0.5em}
 \begin{subfigure}[t]{0.3\textwidth}
  \centering
  \includegraphics[width=\textwidth]{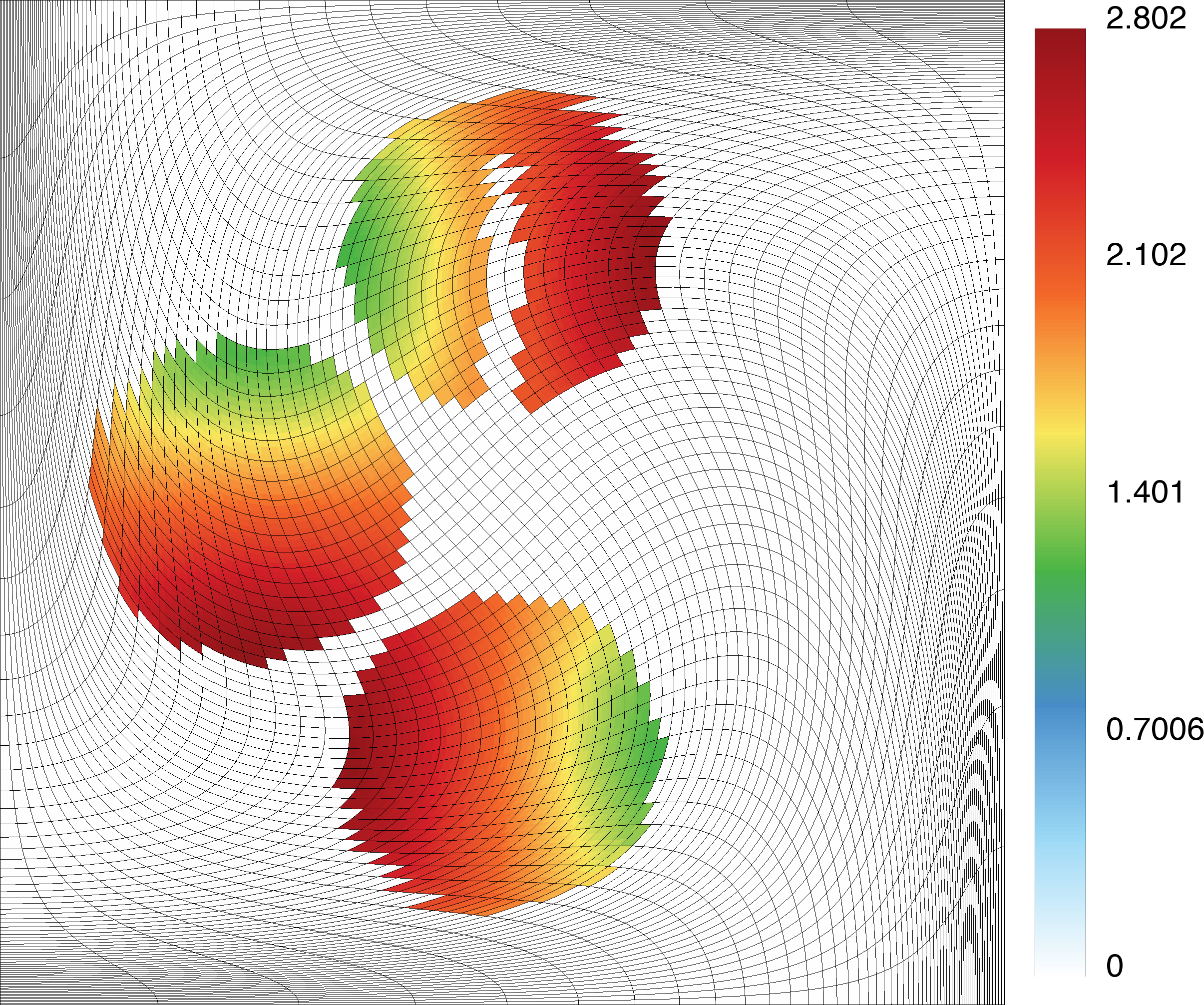}
  \caption{$\rho$, remapped}
\end{subfigure}\hfil
 \begin{subfigure}[t]{0.3\textwidth}
  \centering
  \includegraphics[width=\textwidth]{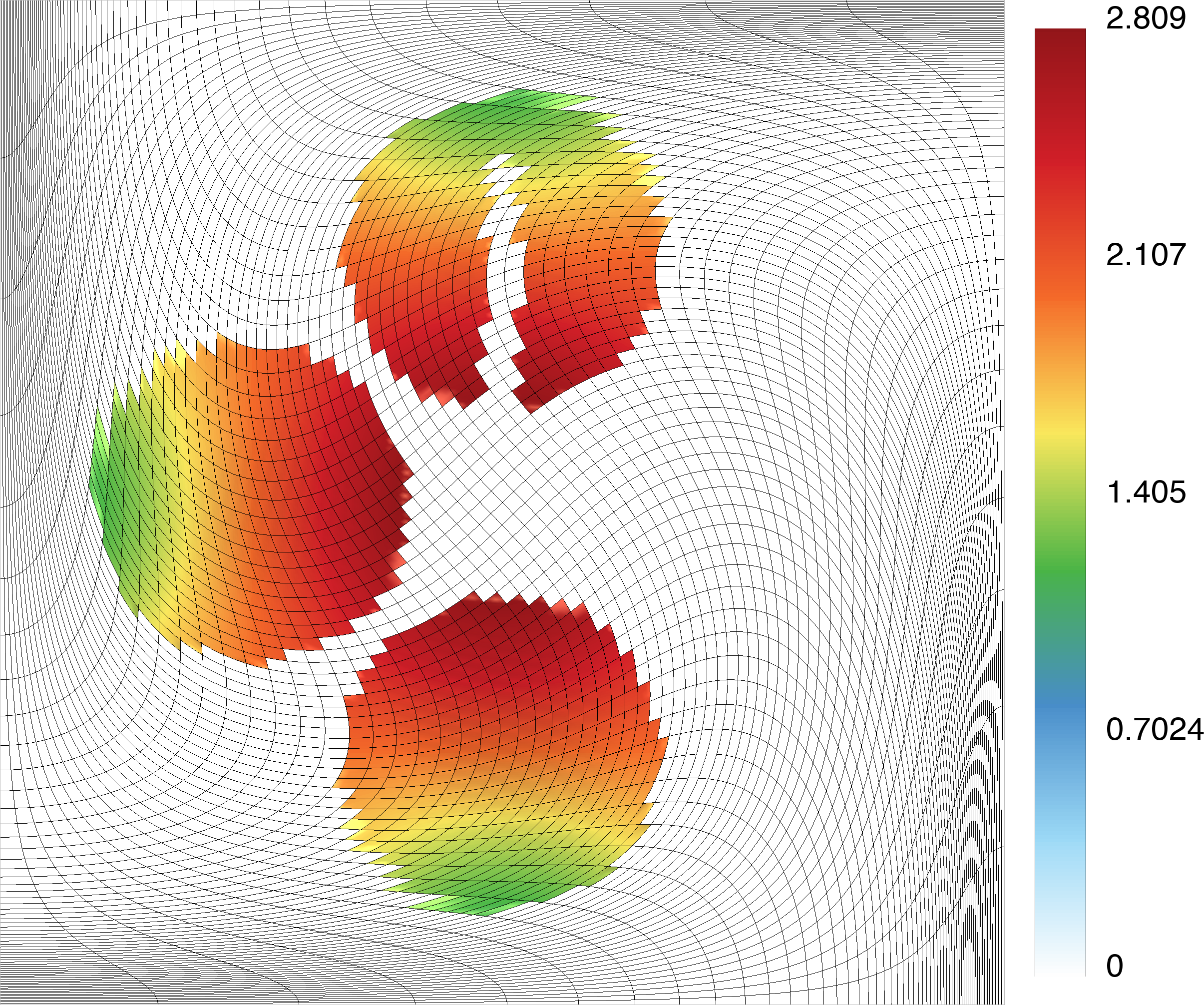}
  \caption{$e$, remapped}
\end{subfigure}\hfil
\begin{subfigure}[t]{0.3\textwidth}
  \centering
  \includegraphics[width=\textwidth]{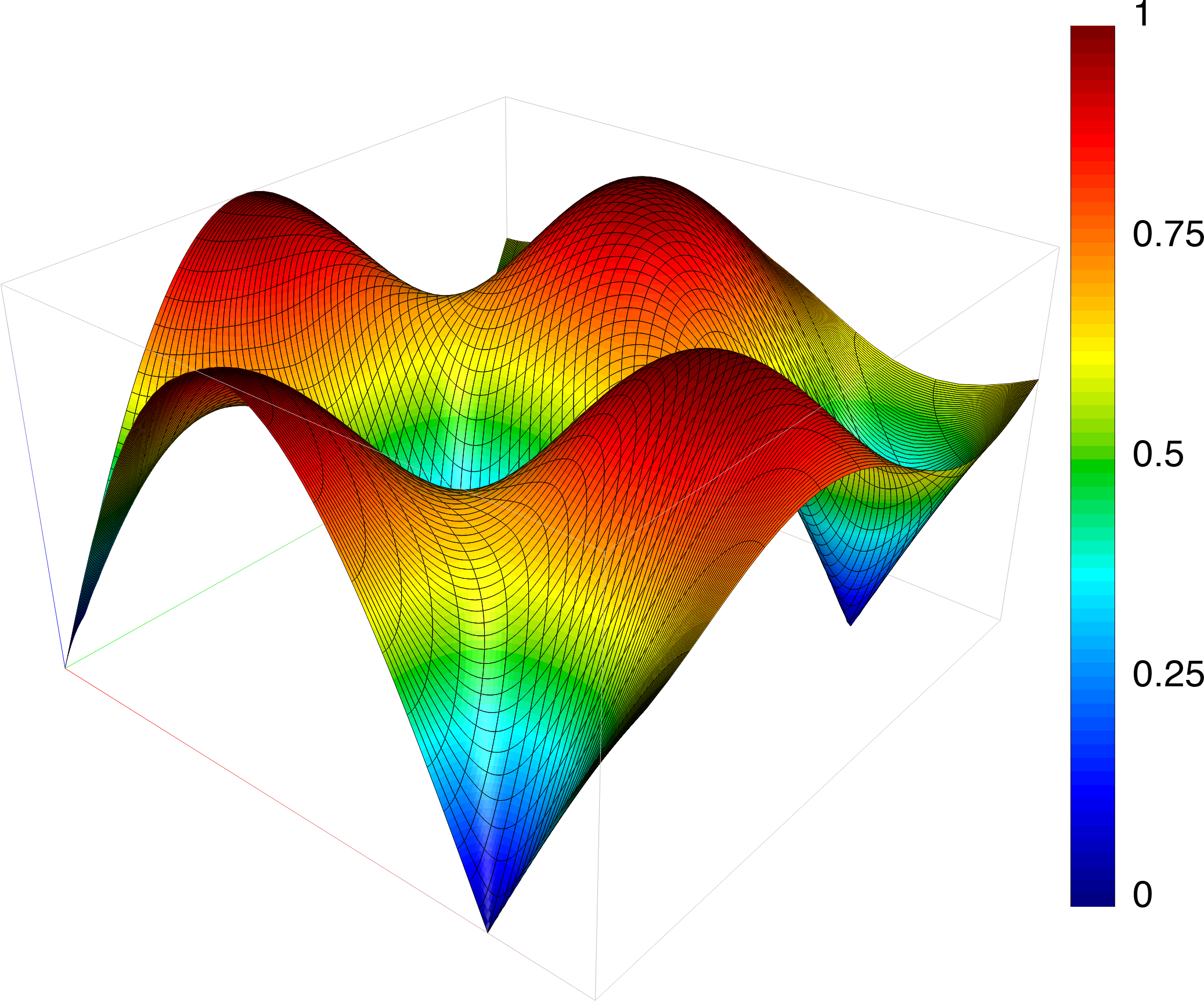}
  \caption{$v$, remapped}
\end{subfigure}
\caption{Initial (top row) and remapped (bottom row) solutions for
         density $\rho$, internal energy $e$, and velocity $v$
         in the Three-Shape benchmark.}
\label{fig_3shape_rev}
\end{figure}

The remaining results demonstrate the coupled-field behavior.
The initial and remapped density, the recovered specific internal energy, and
the velocity fields are presented in Figure~\ref{fig_3shape_rev}.
Their support remains aligned with the material indicator, demonstrating
that the pressure-controlled optimization and subsequent energy recovery
preserve consistent hydrodynamic states near the interface
while preserving the indicator, density, and pressure bounds.
The enforced coupled conservations of mass, momentum, and
total energy are reported in Table~\ref{tab_3shape_conservation}.
The high-order interpolation introduces sub-percent changes in these integral
quantities, whereas the optimized remap reduces the conservation errors to the
prescribed optimization tolerance.

\begin{table}[h!]
\centering
\caption{Conservation errors before and after the optimization step
         for the Three-Shape benchmark.}
\label{tab_3shape_conservation}
\small
  \begin{tabular}{lrrrr}
  \hline
  & \multicolumn{2}{c}{Interpolated} & \multicolumn{2}{c}{Optimized} \\
  Quantity & Abs. error & Rel. error (\%) & Abs. error & Rel. error (\%) \\
  \hline
  Volume         & $+3.60$E$-4$ & $+3.94$E$-1$ & $-4.16$E$-17$ & $-4.56$E$-14$ \\
  Mass           & $+7.92$E$-4$ & $+4.34$E$-1$ & $-5.28$E$-11$ & $-2.89$E$-8$  \\
  Momentum $x$   & $+4.05$E$-4$ & $-5.60$E$-1$ & $-4.28$E$-10$ & $+5.92$E$-7$  \\
  Momentum $y$   & $+6.22$E$-5$ & $+6.20$E$-1$ & $+2.82$E$-11$ & $+2.81$E$-7$  \\
  Total energy   & $+1.84$E$-3$ & $+4.44$E$-1$ & $+8.08$E$-10$ & $+1.95$E$-7$  \\
  \hline
  \end{tabular}
\end{table}

\subsection{Triple Point Problem}
\label{sec_result_triple_point}

The final example considers a single-material variant of the Triple-Point
problem, a standard hydrodynamics benchmark exhibiting interacting shocks and
strong shear that generate complex vortical structures.
The configuration is initialized with a triple-point-style density
discontinuity and uses the same ideal-gas equation of state,
with $\gamma=1.5$, throughout the domain.
This example explores the remap method’s ability to preserve solution quality
in regions dominated by shocks and vorticity over many remap cycles.

The analysis uses a $Q_2Q_1$ discretization and advances to $T=5.0$,
initializing an ALE remap step every $0.5$ time units.
A key ingredient for this test is pressure control during remap:
pressure is remapped as a bounded Q-function at quadrature points,
and the specific internal energy is recovered from the remapped density
and pressure through the equation of state.
This keeps pressure within bounds during remap, which makes
a significant difference in the robustness of the overall ALE evolution and
allows the simulation to complete without developing nonphysical instabilities.
In addition, the numerical diffusion introduced by the remap remains confined
to a narrow band of approximately one to two elements across the density jump,
consistent with the intended sharpness of the method.
The results shown here use the Bregman-projection method
from Section \ref{sec_opt_bp}.
A corresponding analysis using the IP method from Section~\ref{sec_opt_hiop}
yields comparable solution quality and interface sharpness.

\begin{figure}[pos=htbp]
\centering
\centerline
{
  \includegraphics[height=0.28\textwidth]{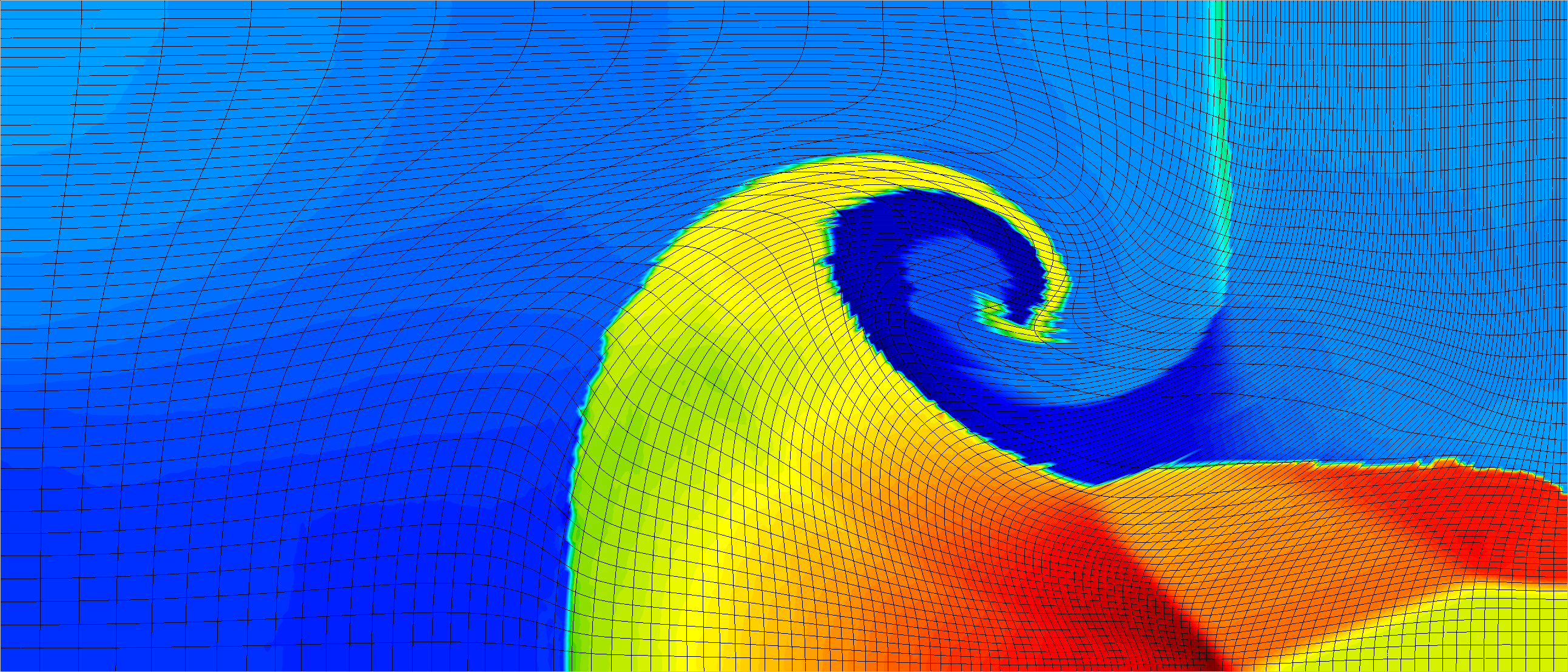} \hfil
  \includegraphics[height=0.28\textwidth]{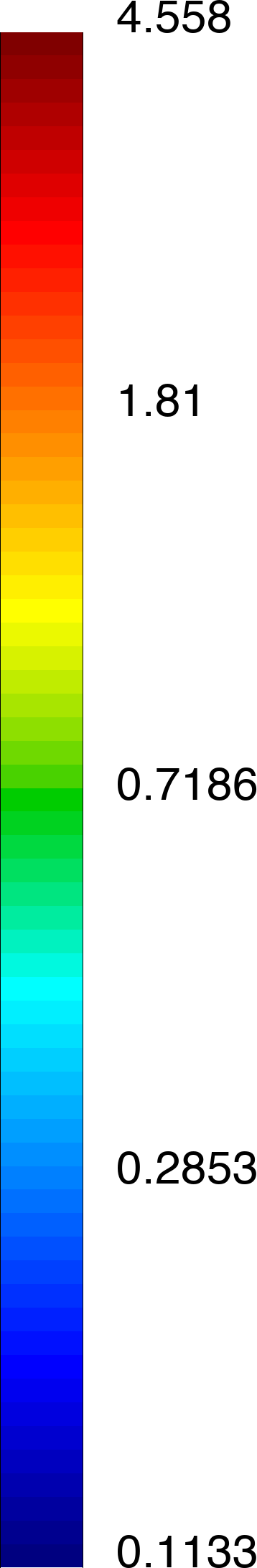}
}
\caption{Density $\rho$ for the two dimensional Triple Point problem
         at $T = 5.0$ using a $Q_2Q_1$ discretization with
         ALE remap every $0.5$ time units (left)
         and the corresponding color scale (right).}
\label{fig_3point_rho}
\end{figure}

\begin{figure}[pos=htbp]
\centering
\includegraphics[width=0.5\textwidth]{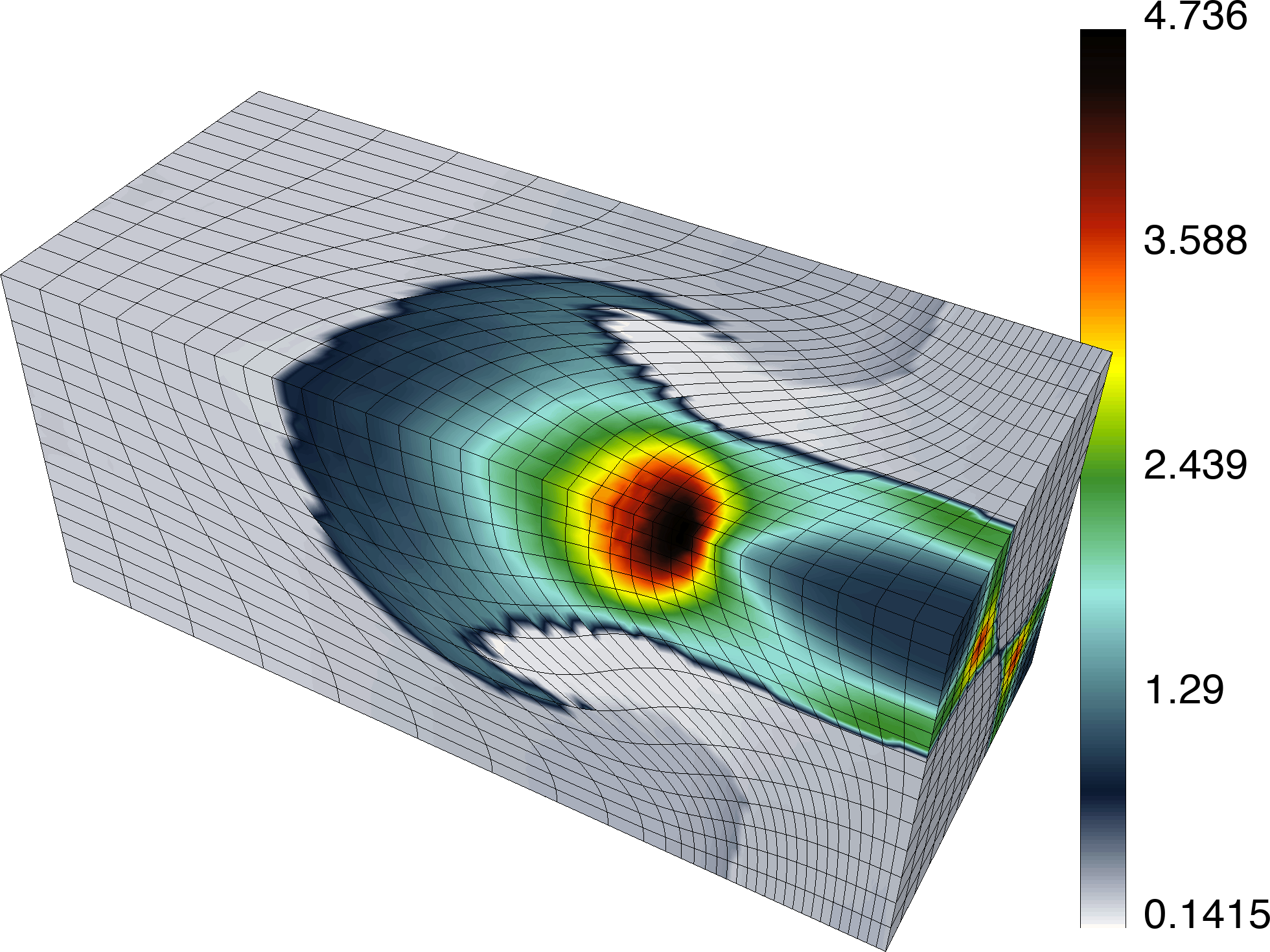}
\caption{3D Triple Point result for the one-material configuration.}
\label{fig_3point_3d}
\end{figure}

The three-dimensional result in Figure~\ref{fig_3point_3d}
demonstrates the method’s direct extension to 3D.
The optimization formulation and its conservation and bound constraints
operate on discrete field vectors, allowing the same algorithmic components
to be used without dimension-specific modifications.
The simulation advances to $T=5.0$,
with an ALE remap step initialized every $1.0$ time units.
Interface sharpness is comparable to the two dimensional results,
with diffusion confined to a band region of approximately one to two elements.

\section{Conclusion}
\label{sec_concl}
This work presents an optimization-based remapping approach for single-material
ALE hydrodynamics that simultaneously enforces global conservation constraints
while preserving local bounds on the
material indicator, density, and auxiliary pressure fields.
The approach combines direct interpolation in physical space with an
optimization-based correction to enforce global conservation of
material volume, mass, momentum, and total energy.
The method maintains sharp interfaces with minimal diffusion and provides a
robust treatment of coupled quantities across different discretizations:
Q-functions for the indicator, density, pressure, and internal energy, and a
high-order $H^1$ vector FE function for velocity.

The formulation offers greater flexibility than traditional advection-based
approaches in controlling pressure and small indicator values, which is critical
for robustness in shock-dominated ALE simulations with sharp transitions.
The numerical results demonstrate that this remap method is
conservative without compromising high-order convergence.
The benchmark tests confirm that the
coupled fields remain aligned with the material indicator and
that optimization reduces conservation errors to the prescribed tolerance.
The two- and three-dimensional triple-point simulations further
demonstrate robustness over repeated remap cycles.

Future work will extend the method to
fully multi-material settings by enforcing inter-material consistency.
This includes partition-of-unity constraints for multiple indicator fields,
along with the conservation and boundedness properties of the present method.
Further directions include improving the treatment of strong discontinuities
over repeated remap cycles and extending high-order,
property-preserving remap to additional point-based quantities.

\section*{Acknowledgements}
\noindent
This manuscript has been authored by Lawrence Livermore National Security,
LLC under Contract No. DE-AC52-07NA27344 with the US. Department of Energy
and was supported by the LLNL-LDRD Program under project No 25-ERD-030.
The United States Government retains, and the publisher, by accepting the
article for publication, acknowledges that the United States Government
retains a non-exclusive, paid-up, irrevocable, world-wide license to publish
or reproduce the published form of this manuscript, or allow others to do so,
for United States Government purposes. LLNL-JRNL-2018243.

BK and DK were supported in part by the U.S. Department
of Energy, Office of Science Early Career Research Program under Award Number
DE-SC0024335 and by the Center for Information Geometric Mechanics and
Optimization (CIGMO), a PSAAP-IV Focused Investigatory Center funded by the
U.S. Department of Energy, National Nuclear Security Administration under Award
Number DE-NA0004261. BK was also supported in part by the Alfred P. Sloan
Foundation via a Sloan Research Fellowship in Mathematics.

\bibliographystyle{cas-model2-names}
\bibliography{refs}

@Article{Svanberg2002,
  Title                    = {A Class of Globally Convergent Optimization Methods Based on Conservative Convex Separable Approximations},
  Author                   = {Krister Svanberg},
  Journal                  = {SIAM Journal on Optimization},
  Year                     = {2002},
  Pages                    = {555-573},
  Volume                   = {12}
}

@article{Levy2024,
author = {Levý, Tomáš and May, Georg},
title = {Conservative solution transfer between anisotropic meshes for time-accurate hybridized discontinuous Galerkin methods},
journal = {International Journal for Numerical Methods in Fluids},
volume = {96},
number = {6},
pages = {1011-1030},
doi = {https://doi.org/10.1002/fld.5278},
url = {https://onlinelibrary.wiley.com/doi/abs/10.1002/fld.5278},
eprint = {https://onlinelibrary.wiley.com/doi/pdf/10.1002/fld.5278},
year = {2024}
}

@article{Bussetta2015,
author = {Bussetta, Philippe and Boman, Romain and Ponthot, Jean-Philippe},
title = {Efficient 3D data transfer operators based on numerical integration},
journal = {International Journal for Numerical Methods in Engineering},
volume = {102},
number = {3-4},
pages = {892-929},
doi = {https://doi.org/10.1002/nme.4821},
url = {https://onlinelibrary.wiley.com/doi/abs/10.1002/nme.4821},
eprint = {https://onlinelibrary.wiley.com/doi/pdf/10.1002/nme.4821},
year = {2015}
}

@article{Scrimieri2014,
title = {Fast mapping of finite element field variables between meshes with different densities and element types},
journal = {Advances in Engineering Software},
volume = {67},
pages = {90-98},
year = {2014},
issn = {0965-9978},
doi = {https://doi.org/10.1016/j.advengsoft.2013.08.003},
url = {https://www.sciencedirect.com/science/article/pii/S096599781300152X},
author = {Daniele Scrimieri and Shukri M. Afazov and Adib A. Becker and Svetan M. Ratchev}
}

@inproceedings{Remap2026,
author    = "Brendan Keith and Dohyun Kim and Boyan Lazarov and
             Cosmin Petra and Mathias Schmidt and Vladimir Z. Tomov",
title     = "Remapping through direct interpolation and optimization for
             finite element {ALE} hydrodynamics",
booktitle = "Large-Scale Scientific Computations",
pages     = "13-21",
year      = "2026",
publisher = "Springer Nature Switzerland",
editor    = "Lirkov, Ivan and Margenov, Svetozar",
doi       = "10.1007/978-3-032-22221-3\_2"
}

@article{MFEM2024,
author  = "Julian Andrej and Nabil Atallah and Jan-Phillip B{\"a}cker and
           John Camier and Dylan Copeland and Veselin Dobrev and
           Yohann Dudouit and Tobias Duswald and Brendan Keith and
           Dohyun Kim and Tzanio Kolev and Boyan Lazarov and
           Ketan Mittal and Will Pazner and Socratis Petrides and
           Syun'ichi Shiraiwa and Mark Stowell and Vladimir Tomov",
title   = "High-performance finite elements with {MFEM}",
journal = "Int. J. High Perform. Comput. Appl.",
volume  = "38",
number  = "5",
pages   = "447-467",
year    = "2024",
doi     = "10.1177/10943420241261981"
}

@article{Kuzmin2024,
author  = "Falko Ruppenthal and Denis Ridzal and Dmitri Kuzmin and Pavel Bochev",
title   = "Scalable optimal control for inequality-constrained
           discretizations of scalar conservation laws",
journal = "Math. Control Relat. Fields",
pages   = "",
year    = "2024",
doi     = "10.3934/mcrf.2024068"
}

@article{keith2026proximal,
  title={Proximal {G}alerkin: {A} structure-preserving finite element method for pointwise bound constraints},
  author={Keith, Brendan and Surowiec, Thomas M},
  journal={Foundations of Computational Mathematics},
  volume={26},
  number={1},
  pages={385--481},
  year={2026},
  publisher={Springer}
}

@article{dokken2025latent,
  title={The latent variable proximal point algorithm for variational problems with inequality constraints},
  author={Dokken, J{\o}rgen S and Farrell, Patrick E and Keith, Brendan and Papadopoulos, Ioannis PA and Surowiec, Thomas M},
  journal={Computer Methods in Applied Mechanics and Engineering},
  volume={445},
  pages={118181},
  year={2025},
  publisher={Elsevier}
}

@article{Bochev2020,
author = "Pavel Bochev and Denis Ridzal and Marta D'Elia and
          Mauro Perego and Kara Peterson",
title  = "Optimization-based, property-preserving finite element methods
          for scalar advection equations and their connection to
          {A}lgebraic {F}lux {C}orrection",
journal = "Comput. Methods Appl. Mech. Eng.",
volume  = "367",
pages   = "112982",
year    = "2020",
doi     = "10.1016/j.cma.2020.112982"
}

@Article{Petra2019,
author  = "Cosmin G. Petra",
title   = "A memory-distributed quasi-{N}ewton solver for
           problems with a small number of general constraints",
journal = "J. Parallel and Distrib. Comput.",
volume  = "133",
pages   = "337-348",
year    = "2019"
}

@article{Tomov2018,
author  = "Robert W. Anderson and Veselin A. Dobrev and Tzanio V. Kolev and
           Robert N. Rieben and Vladimir Z. Tomov",
title   = "High-Order Multi-Material {ALE} Hydrodynamics",
journal = "SIAM J. Sci. Comp.",
volume  = "40",
number  = "1",
pages   = "B32-B58",
year    = "2018",
doi     = "10.1137/17M1116453"
}

@article{Manuel2017,
author  = "R. W. Anderson and V. A. Dobrev and T. V. Kolev and D. Kuzmin and
           M. {Quezada de Luna} and R. N. Rieben and V. Z. Tomov",
title   = "High-Order Local Maximum Principle Preserving ({MPP}) Discontinuous
           {G}alerkin Finite Element Method for the Transport Equation",
journal = "J. Comput. Phys.",
volume  = "334",
pages   = "102-124",
year    = "2017",
doi     = "10.1016/j.jcp.2016.12.031"
}

@article{Dobrev2015,
author  = "R. W. Anderson and V. A. Dobrev and T. V. Kolev and R. N. Rieben",
title   = "Monotonicity in high-order curvilinear finite element arbitrary
           {L}agrangian-{E}ulerian remap",
journal = "Int. J. Numer. Methods Fluids",
volume  = "77",
number  = "5",
pages   = "249-273",
year    = "2015",
doi     = "10.1002/fld.3965"
}

@article{Bochev2014,
author = "Pavel Bochev and Denis Ridzal and Kara Peterson",
title  = "Optimization-based remap and transport: a divide and conquer
          strategy for feature-preserving discretizations",
journal = "J. Comput. Phys.",
volume  = "257",
pages   = "1113-1139",
year    = "2014",
doi     = "10.1016/j.jcp.2013.03.057"
}

@article{Bochev2013,
author = "Pavel Bochev and Denis Ridzal and Mikhail Shashkov",
title  = "Fast optimization-based conservative remap of scalar fields
          through aggregate mass transfer",
journal = "J. Comput. Phys.",
volume  = "246",
pages   = "37-57",
year    = "2013",
doi     = "10.1016/j.jcp.2013.03.040"
}

@article{Dobrev2012,
author  = "V. Dobrev and Tz. Kolev and R. Rieben",
title   = "High-Order Curvilinear Finite Element Methods
           for {L}agrangian Hydrodynamics",
journal = "SIAM J. Sci. Comp.",
volume  = "34",
number  = "5",
pages   = "606-641",
year    = "2012",
doi     = "10.1137/120864672"
}

@article{Bochev2011,
author = "Pavel Bochev and Denis Ridzal and
          Guglielmo Scovazzi and Mikhail Shashkov",
title  = "Formulation, analysis and numerical study of an optimization-based
          conservative interpolation (remap) of scalar fields for an
          arbitrary {L}agrangian-{E}ulerian methods",
journal = "J. Comput. Phys.",
volume  = "230",
pages   = "5199-5225",
year    = "2011",
doi     = "10.1016/j.jcp.2011.03.017"
}

@article{Liska2010,
author = "Richard Liska and Mikhail Shashkov and
          Pavel V\'{a}chal and Burton Wendroff",
title  = "Optimization-based synchronized flux-corrected conservative
          interpolation (remapping) of mass and momentum for
          arbitrary {L}agrangian-{E}ulerian methods",
journal = "J. Comput. Phys.",
volume  = "229",
pages   = "1467-1497",
year    = "2010",
doi     = "10.1016/j.jcp.2009.10.039"
}

@article{LeVeque1996,
author = "Randall J LeVeque",
title  = "High-resolution conservative algorithms for advection in
          incompressible flow",
journal = "SIAM J. Numer. Anal.",
volume  = "33",
pages   = "627-665",
year    = "1996",
doi     = "10.1137/0733033"
}

@misc{gslib,
key    = "GSLIB",
author = "Paul Fischer",
title  = "{GSLIB}: sparse communication library {[Software]}",
year   = "2017",
note   = "{\url{https://github.com/gslib/gslib}}"
}

@article{Mittal2025,
  title={General field evaluation in high-order meshes on GPUs},
  author={Mittal, Ketan and Parik, Aditya and Dutta, Som and Fischer, Paul and Kolev, Tzanio and Lottes, James},
  journal={Computers \& Fluids},
  volume={295},
  pages={106606},
  year={2025},
  publisher={Elsevier}
}

@article{waecther_05_ipopt,
author = {Andreas W{\"a}chter and Lorenz T. Biegler},
title = {Line Search Filter Methods for Nonlinear Programming: Motivation and Global Convergence},
journal = {SIAM Journal on Optimization},
volume = {16},
number = {1},
pages = {1-31},
year = {2005}
}

@article{waecther_05_ipopt0,
year={2006},
issn={0025-5610},
journal={Mathematical Programming},
volume={106},
number={1},
title={On the implementation of an interior-point filter line-search algorithm for nonlinear programming},
publisher={Springer-Verlag},
author={W{\"a}chter, Andreas and Biegler, Lorenz T.},
pages={25-57},
language={English}
}

@article{Bregman1967,
  title = {The relaxation method of finding the common point of convex sets and its application to the solution of problems in convex programming},
  volume = {7},
  ISSN = {0041-5553},
  DOI = {10.1016/0041-5553(67)90040-7},
  number = {3},
  journal = {USSR Computational Mathematics and Mathematical Physics},
  publisher = {Elsevier BV},
  author = {Bregman,  L.M.},
  year = {1967},
  month = jan,
  pages = {200–217}
}

@article{Dowell1971,
  title = {A modified regula falsi method for computing the root of an equation},
  volume = {11},
  ISSN = {1572-9125},
  url = {http://dx.doi.org/10.1007/BF01934364},
  DOI = {10.1007/bf01934364},
  number = {2},
  journal = {BIT},
  publisher = {Springer Science and Business Media LLC},
  author = {Dowell,  M. and Jarratt,  P.},
  year = {1971},
  month = jun,
  pages = {168–174}
}

@article{kim2025simple,
  title={A simple introduction to the {SiMPL} method for density-based topology optimization},
  author={Kim, Dohyun and Lazarov, Boyan S and Surowiec, Thomas M and Keith, Brendan},
  journal={Structural and Multidisciplinary Optimization},
  volume={68},
  number={4},
  pages={74},
  year={2025},
  publisher={Springer}
}

@misc{HiOpRepo,
  author       = {{Lawrence Livermore National Laboratory}},
  title        = {HiOp: An optimization solver for large-scale nonlinear programming},
  howpublished = {\url{https://github.com/LLNL/hiop}},
  note         = {Accessed April 3, 2026}
}

@article{Bauschke2000,
  title = {Dykstras algorithm with {B}regman projections: A convergence proof},
  volume = {48},
  ISSN = {1029-4945},
  url = {http://dx.doi.org/10.1080/02331930008844513},
  DOI = {10.1080/02331930008844513},
  number = {4},
  journal = {Optimization},
  publisher = {Informa UK Limited},
  author = {Bauschke,  Heinz H and Lewis,  Adrian S},
  year = {2000},
  month = Jan,
  pages = {409–427}
}

@techreport{keith2025remap,
  title       = {Remapping Through Direct Interpolation and Optimization for Finite Element {ALE} Hydrodynamics},
  author      = {Brendan Keith and Dohyun Kim and Boyan Lazarov and Cosmin Petra and Mathias Schmidt and Vladimir Z Tomov},
  institution = {Lawrence Livermore National Laboratory},
  year        = {2025},
  number      = {LLNL-CONF-2003056},
  url         = {https://www.osti.gov/servlets/purl/2585170},
  address     = {Livermore, CA},
}

@misc{Laghos,
key   = "Laghos",
title = "{L}aghos: High-Order {L}agrangian Hydrodynamics Miniapp",
year  = "2026",
note  = "{https://github.com/CEED/Laghos}"
}

@misc{Remhos,
key   = "Remhos",
title = "{R}emhos: {R}emap High-Order Solver",
year  = "2026",
note  = "{\url{http://github.com/CEED/Remhos}}"
}

@article{anderson,
author = {Anderson, Donald G.},
title = {Iterative Procedures for Nonlinear Integral Equations},
year = {1965},
issue_date = {Oct. 1965},
publisher = {Association for Computing Machinery},
address = {New York, NY, USA},
volume = {12},
number = {4},
issn = {0004-5411},
url = {https://doi.org/10.1145/321296.321305},
doi = {10.1145/321296.321305},
journal = {J. ACM},
month = oct,
pages = {547–560},
numpages = {14}
}

@article{TMOP2019,
author  = "Veselin A. Dobrev and Patrick Knupp and Tzanio V. Kolev and
           Ketan Mittal and Vladimir Z. Tomov",
title   = "The {T}arget-{M}atrix {O}ptimization {P}aradigm
           for High-Order Meshes",
journal = "SIAM J. Sci. Comp.",
volume  = "41",
number  = "1",
pages   = "B50-B68",
year    = "2019",
doi     = "10.1137/18M1167206"
}

@article{MFEM2020,
author  = "Robert Anderson and Julian Andrej and Andrew Barker and
           Jamie Bramwell and Jean-Sylvain Camier and Jakub Cerveny and
           Veselin A. Dobrev and Yohann Dudouit and Aaron Fisher and
           Tzanio V. Kolev and Will Pazner and Mark Stowell and
           Vladimir Z. Tomov and Ido Akkerman and Johann Dahm and
           David Medina and Stefano Zampini",
title   = "{MFEM}: a modular finite elements methods library",
journal = "{Computers \& Mathematics with Applications}",
volume  = "81",
pages   = "42-74",
year    = "2020",
doi     = "10.1016/j.camwa.2020.06.009"
}

\end{document}